\documentclass[12pt]{amsart}
\usepackage{amsmath}
\usepackage{amsfonts}
\usepackage{amssymb}
\usepackage{graphicx}
\usepackage{color}

\numberwithin{equation}{section}
\definecolor{mycolorred}{rgb}{1, 0, 0}
\def\red #1{{#1}}
\definecolor{mycolorblue}{rgb}{0, 0, 1}

\definecolor{mycolorpink}{rgb}{1, 0, 1}

\newtheorem{theorem}{Theorem}[section]

\newtheorem{definition}[theorem]{Definition}

\newtheorem{lemma}[theorem]{Lemma}

\newtheorem{proposition}[theorem]{Proposition}
\newtheorem{remark}[theorem]{Remark}

\def\<{\langle}
\def\>{\rangle}

\def\P{{\mathbb P}}

\def\E{{\mathbb E}}
\def\R{{\mathbb R}}

\def\N{{\mathbb N}}
\def\Z{{\mathbb Z}}

\def\DD{\red{\mathcal{D}}}

\begin{document}
\def\kp{{\kappa_p}}

\title{Stochastic sewing lemma on Wasserstein space}
\author{Aur\'elien Alfonsi}
\address{CERMICS, Ecole des Ponts, Marne-la-Vall\'ee, France. MathRisk, Inria, \
Paris,
  France.}
\email{aurelien.alfonsi@enpc.fr}
\author{Vlad Bally}
\address{Universit\'e Gustave Eiffel, LAMA (UMR CNRS, UPEMLV, UPEC), MathRisk INRIA,
  F-77454 Marne-la-Vall\'ee, France.}
\email{vlad.bally@univ-eiffel.fr}
\author{Lucia Caramellino}
\address{Dipartimento di matematica, Universit\`a degli studi di Roma Tor Vergata -- Via della ricerca scientifica 1, 00133 Roma}
\email{caramell@mat.uniroma2.it}

\begin{abstract}
  The stochastic sewing lemma recently introduced by Le~\cite{[Le]} allows to construct a unique limit process from a doubly indexed stochastic process that satisfies some regularity. This lemma is stated in a given probability space on which these processes are defined. The present paper develops a version of this lemma for probability measures: from a doubly indexed family of maps on the set of probability measures that have a suitable probabilistic representation, we are able to construct a limit flow of maps on the probability measures. This result complements and improves the existing result coming from the classical sewing lemma. It is applied to the case of law-dependent jump SDEs for which we obtain weak existence result as well as the uniqueness of the marginal laws. 
\end{abstract}
\maketitle

\noindent {\bf Acknowledgment:} { This research is partly funded by the Bézout Labex, funded by ANR, reference ANR-10-LABX-58. A.A. acknowledges the support of the "chaire Risques Financiers", Fondation du Risque. L.C. acknowledges the MUR Excellence Project MatMod@TOV awarded to the Department of Mathematics of the University of Rome Tor Vergata.
} \\
\parindent 0pt

\bigskip

\section{Introduction}

The sewing lemma has been introduced in the early 2000 by Gubinelli~\cite{[G]} and
by De la Pradelle and Feyel~\cite{[FP]}. It has 
shed a new insight  on the theory of rough paths, see in particular Davie~\cite{[D]} and the monograph by Friz and Hairer~\cite{[FH]}. A long list of variants and generalizations have  followed, in particular the recent ones of Brault and Lejay~\cite{BrLe1,BrLe2,BrLe3} and it is still nowadays an extremely active research area. In
particular, we mention some variants of the sewing lemma in connection with
flows: a first result in his direction is due to Bailleul~\cite{[B]}, and  more recently Alfonsi and Bally~\cite{[AB]} have used the sewing lemma to construct flows on the probability measures, in
order to construct after solutions for non linear SDEs of McKean and
Boltzmann type. The present paper goes in the same direction, but relies on the new results of L\^e~\cite{[Le]} on the so-called stochastic sewing lemma.

Let us stress a specific aspect of this theory. A crucial point in the rough
path theory is to obtain a pathwise construction of the stochastic
integral which avoids any probabilistic argument. Let us be more explicit. The
stochastic integral (and more generally, martingales) moves with $\sqrt{h}$ in
a time step $h$. Then, in order to avoid the blow up of Riemann sums that define the stochastic integral, one has to
use the fact that the increments are centred and so orthogonal in $L^{2}$.
This is the probabilistic argument which allows to construct the stochastic
integral as a limit in $L^{2}$. In contrast, in order to obtain the pathwise
construction, Lyons introduced a correction given by the L\'evy area, which
allows to handle the corrected Riemann sums. It is precisely in the analysis of 
 these sums that the (classical) sewing lemma is used.

Having this in mind, it is rather striking that in the stochastic sewing
lemma developed by L\^e~\cite{[Le]} comes back and introduces once again the probabilistic argument
based on the martingale property. The fact is that, in contrast with the
previous approach related to rough paths, a probabilistic structure given by a
filtered probability space $(\Omega, \mathcal{F}, \P, (\mathcal{F}_t)_{t\ge 0})$ and discrete time martingales are underlying this
procedure. The basic object in the stochastic sewing lemma setting is a family of $d$ dimensional, square integrable random variables $A_{s,t},s<t,$ which are
adapted to the above mentioned filtration and which, roughly speaking,
represent an one step approximation scheme from $s$ to $t$ that can be typically the Euler scheme. The ``sewing error"
that controls the convergence of the scheme is $\delta
 A_{s,r,t}=A_{s,t}-(A_{s,r}+A_{r,t}), \ s<r<t$, and~\cite{[Le]}  distinguishes a martingale part $\delta A_{s,r,t}- \E[\delta A_{s,r,t}|\mathcal{F}_s]$ and a ``finite variation" part $\E[\delta A_{s,r,t}|\mathcal{F}_s]$
that is analogous to the Doob-Meyer decomposition. Suitable estimates are assumed on these two parts, and the martingale part plays the crucial role in the convergence argument to get the limit object.

The aim of our paper is to give a variant of the stochastic sewing lemma on
the Wasserstein space, and to use it then to construct flows on the space of probability measures. This is a first step to get existence  and uniqueness results for non linear SDEs, as we show on a generic example in Section~\ref{Sec_jump}. We
consider a Hilbert space $H$ (in particular we may take $H=\R^{d}$) and for $p\ge 2$, we
denote by $\mathcal{P}_{p}(H)$ the space of probability measures $\mu$
on $(H,\mathcal{B}(H))$ which have a finite moment of order~$p$. Then, the
basic objet in our approach is a family of applications $\Theta _{s,t}:%
\mathcal{P}_{p}(H)\rightarrow \mathcal{P}_{p}(H),$ $s<t,$ which represents distributions of a discrete approximation of a dynamics. For example, $\Theta _{s,t}(\mu)$ may be the distribution of the one step Euler scheme of a McKean-Vlasov SDE starting from the initial distribution~$\mu$. 
Here is the specific difference with the work of L\^e~\cite{[Le]}: while~\cite{[Le]} considers a filtered probability space is given with a family of random variables $A_{s,t}$, our setting brings on  probability measures $\mu $ and $%
\Theta _{s,t}(\mu )$, and so we do not have access to
an underlying martingale structure. In order to get around this difficulty we assume that a representation property holds: this is given by the ``$\beta$-coupling family of operators" defined in Definition~\ref{def_betacoupling} for $\beta>0$. Once this representation is settled, we follow the martingale approach (analogous but not identical to the one of~\cite{[Le]}), combined with optimal coupling arguments.
To state our version of the stochastic sewing lemma at the level of the Wasserstein space (i.e. in terms of the family $\Theta _{s,t}$), we consider partitions. More precisely, we define the scheme $\Theta_{s,t}^{\pi}(\mu )=\Theta _{t_{N-1},t_{N}}\circ ...\circ \Theta
_{t_{0},t_{1}}(\mu )$ associated with the partition $\pi
=\{s=t_{0}<...<t_{N}=t\}$. We show in Lemma~\ref{sewing} that $W_p(\Theta^{\pi}(\mu),\Theta^{\pi'}(\mu))=O(\max_{0\le i\le n}|t_i-t_{i-1}|^\beta)$ when $\pi'=\{t_{0}\leq r_{0}<t_{1}<....<t_{N-1}\leq r_{N-1}<t_{N}\}$ is a partition obtained from~$\pi$ by adding at most one point on each interval $[t_i,t_{i+1})$.  In contrast with the classical sewing lemma that basically assumes the local estimate $W_p(\Theta_{s,t}(\mu), \Theta_{r,t}\circ \Theta_{s,r}(\mu))=o(t-s)$, this stochastic sewing lemma quantifies the Wasserstein distance at the level of partitions where the martingale argument is nested. From this lemma, we then prove our 
main result (Theorem~\ref{existence}) and get that $\Theta _{s,t}^{\pi}(\mu
)\rightarrow \theta _{s,t}(\mu )$ as $\left\vert \pi\right\vert \rightarrow
0$, where $\theta _{s,t}:\mathcal{P}_{p}(H)\rightarrow \mathcal{P}_{p}(H)$ is a flow: $\theta _{s,t}=\theta _{r,t}\circ
\theta _{s,r}$ for $s<r<t$. Besides, the stochastic sewing lemma  allows then to give a speed of convergence for $W_p(\Theta_{s,t}^{\pi}(\mu),\theta_{s,t}(\mu))$ in function of the time step size of the partition~$\pi$, see Theorem~\ref{main}.

The paper is structured in two sections. Section~\ref{Sec_sewing_flow} presents the general framework, the stochastic sewing lemma on Wasserstein space and the construction of the flow $\theta_{s,t}$ from the scheme $\Theta_{s,t}$. Section~\ref{Sec_jump} presents then an application of this theorem to non-linear McKean-Vlasov SDEs with jumps. More precisely, our stochastic sewing result is applied to get the marginal laws, which enables us then to get existence and uniqueness for these equations.

\section{Stochastic sewing lemma on the Wasserstein space and flows}\label{Sec_sewing_flow}

Through the paper a separable Hilbert space $H$ is fixed, and we denote by $\left\vert x\right\vert _{H}$ the norm on $H$ and by $\mathcal{B}(H)$ the
Borel $\sigma$-algebra on $H$. We fix $p\geq 2$ and we denote by $\mathcal{P
}_{p}(H)$ the space of probability measures $\mu$ on $(H,\mathcal{B}(H)),$
such that $\int_{H}\left\vert x\right\vert _{H}^{p}\mu (dx)<\infty$. We note in this case
\begin{equation}
\|\mu\|_p=\left(\int_{H}\left\vert x\right\vert _{H}^{p}\mu (dx)\right)^{1/p} \text{ and } \|X\|_p= \E[|X|_H^p]^{1/p}=\|\mu\|_p,
\end{equation}
for any random variable $X$ that has the distribution~$\mu$. We also recall the definition of the Wasserstein distance:
$$ W_p(\mu,\nu)= \inf_{\Pi \in  \text{Cpl}(\mu,\nu)} \left(\int_{H\times H} |x-y|_H^p\Pi(dx,dy)\right)^{1/p},$$
with $\text{Cpl}(\mu,\nu)=\{ \Pi \in \mathcal{P
}_{p}(H\times H): \mu(dx)=\int_{y\in H}\Pi(dx,dy), \nu(dy)=\int_{x\in H}\Pi(dx,dy)\}$.
\begin{remark}\label{Rk_Banach}
  In this paper, we consider \red{a Hilbert} space for convenience. In fact, the arguments developed in this paper would work also in Banach spaces, provided that we have at hand a Burkholder-Davis-Gundy inequality. This point is discussed in \cite[p.~3]{[Le2]}, and a list of Banach spaces for which a BDG inequality holds is given.  
\end{remark}

\begin{definition}\label{def_partition}
  Let $0\le t<u$. A partition $\pi$ of $[t,u]$ is a finite increasing sequence 
  $(t_i)_{0\le i\le N}$ such that $t_0=t$ and $t_N=u$. We note $\pi
  =\{t_{0}<t_{1}<....<t_{N}\}$ and we denote $|\pi|=\max_{i=1,%
  \dots,N}\red{(t_i-t_{i-1})}$ and $\#\pi=N+1$.
  \end{definition}

  We consider a family of applications $\Theta _{s,t}:\mathcal{P}%
  _{p}(H)\rightarrow \mathcal{P}_{p}(H),0\leq s\leq t$ such that $\Theta_{t,t}=
  \mathrm{Id}$. We may think to $\Theta _{s,t}(\mu )$ as a one step Euler
  scheme starting from the initial law~$\mu$ at time $s$. 
  For a partition (often also called time grid) $\pi
=\{t_{0}<t_{1}<....<t_{N}\}$, with $0\leq t_0$, 
we define the operator $\Theta _{t_{0},t_{N}}^{\pi }:\mathcal{P}_{p}(H)\rightarrow \mathcal{P}_{p}(H) $ by
\begin{equation}\label{def_Thetapi}
\Theta_{t_{0},t_{N}}^{\pi }(\mu )=\Theta _{t_{N-1},t_{N}}\circ ...\circ
\Theta_{t_{0},t_{1}}(\mu ).
\end{equation}%
Again, we may think to $\Theta _{t_{0},t_{N}}^{\pi }(\mu )$ as the law of a Euler scheme at time $t_N$ that uses the time grid~$\pi$ and starts from $t_0$ with law~$\mu$. This is the main object of interest in our paper. Our main hypothesis is
that $\Theta _{t_{0},t_{N}}^{\pi }(\mu )$ admits a convenient probabilistic
representation which allows us to separate a ``martingale part" and a ``finite
variation part". We discuss this in the following section. The aim of the
paper is to construct a flow $\theta _{s,t}:\mathcal{P}_{p}(H)\rightarrow 
\mathcal{P}_{p}(H),0\leq s\leq t$ as limit of $\Theta _{s,t}^{\pi }$ as the
mesh $\left\vert \pi \right\vert=\max_{1\le i\le N}|t_i-t_{i-1}|$ of the partition goes to zero.

\subsection{Coupling operator and stochastic sewing lemma}

In this section, we still consider a family of operators $\Theta _{s,t}:
\mathcal{P}_{p}(H)\rightarrow \mathcal{P}_{p}(H)$ and we want to define a \textbf{coupling operator } that gives a joint probabilistic representation of the laws $\Theta_{s,t}(\mu)$ and $\Theta_{r,t}\circ\Theta_{s,r}(\mu)$, for $s<r<t$. To define such an operator, we will use a probability space $(\Omega ,\mathcal{F},\P)$ with a filtration $(\mathcal{F}_{t},t\geq 0)$, which is fixed through this section. We assume that 
\begin{itemize}
  \item[\textbf{A}] \textbf{[Atomless]} There   exists a $\mathcal{F}_0$-measurable random variable $U:\Omega \to \R$ such that $U$ has the uniform distribution on $[0,1]$.  This is known to be equivalent to the atomless property of the probability space $(\Omega,\mathcal{F}_0,\P)$.
\end{itemize}  
Then\footnote{By Kuratowski's theorem, every Polish space is Borel isomorphic to either $\R$, $\Z$ or a finite set. }, for every $s\geq 0$ and every
probability measure $\eta$ on $H$ (or on $H\times H$), there exists an $\mathcal{F}_{s}$ (even $\mathcal{F}_0$) measurable random variable $Z$ such that $Z\sim \eta$.

For $t\ge 0$, we denote by $L^p_t(H)$ the space of $\mathcal{F}_t$-measurable random variables~$X:\Omega \to H$ such that $\|X\|_p<\infty$. The main object of this section is a family of operators:
$$\mathcal{X}_{s,r,t}: L_{s}^{p}(H)\times L_{s}^{p}(H)\to  L_{t}^{p}(H)\times L_{t}^{p}(H),
\quad \red{s\leq r\leq t}.$$
For $X=(X^1,X^2)\in L_{s}^{p}(H)\times L_{s}^{p}(H)$, we write 
$$\mathcal{X}_{s,r,t}(X^1,X^2)=\left(\mathcal{X}^1_{s,r,t}(X^1,X^2) ,\mathcal{X}^2_{s,r,t}(X^1,X^2) \right) \in  L_{t}^{p}(H)\times L_{t}^{p}(H).$$
We stress that all these operators use the same filtered probability space, for every $s\leq r \leq t$. 
Let us stress also that $\mathcal{X}_{s,r,t}$  takes in input random variables  (i.e the map $\omega \mapsto (X^1(\omega),X^2(\omega))$), not two elements of~$H$. This gives in particular the possibility to use the joint law $\mathcal{L}(X^1,X^2)$ in the definition  of $\mathcal{X}_{s,r,t}$ since it is the image measure (also called pushforward measure) of~$\mathbb{P}$ (that is fixed) by the map $(X_1,X_2)$. 

Let us denote, for $\ell \in \{1,2\}$, $\mathcal{E}_{s,r,t}^{\ell}(X)=\mathcal{X}_{s,r,t}^{\ell}(X)-\red{X^\ell}$ and 
\begin{equation}
\widehat{\mathcal{E}}_{s,r,t}^{\ell}(X)=\E_{\mathcal{F}_{s}}(\mathcal{E}%
_{s,r,t}^{\ell}(X)),\quad \widetilde{\mathcal{E}}_{s,r,t}^{\ell}(X)=\mathcal{E}%
_{s,r,t}^{\ell}(X)-\widehat{\mathcal{E}}_{s,r,t}^{\ell}(X),  \label{def_E_hat_tilde}
\end{equation}
where $\E_{\mathcal{F}_{s}}$ denotes the conditional expectation with respect to~$\mathcal{F}_s$.

We assume that these operators verify the following properties.

\begin{itemize}
\item[\textbf{R}] [\textbf{Representation}]  For every $s\leq r\leq t$,
$\mu ^{1},\mu ^{2}\in \mathcal{P}_{p}(H)$ and every $X^{1},X^{2}\in
L_{s}^{p}(H)$ with $\mathcal{L}(X^{i})=\mu ^{\red{i}},\red{i}=1,2$,%
\begin{equation*}
\Theta _{s,t}(\mu ^{1})=\mathcal{L}(\mathcal{X}_{s,r,t}^{1}(X^{1}, \red{X^2}))\quad and\quad
\Theta _{r,t}\circ \Theta _{s,r}(\mu ^{2})=\mathcal{L}(\mathcal{X}_{s,r,t}^{2}(\red{X^1},X^{2})).
\end{equation*}
\item[$\mathbf{G}$] [\textbf{Growth assumption}] There exists $L_{sew}>0$ such that  \red{for every $s\leq r\leq t$}
\begin{equation}
\left\Vert \widehat{\mathcal{E}}^1_{s,r,t}(X)\right\Vert _{p}\leq
L_{sew}(t-s)(1+\left\Vert X^1\right\Vert _{p})\quad \left\Vert \widetilde{%
\mathcal{E}}_{s,r,t}^{1}(X)\right\Vert _{p}\leq
L_{sew}(t-s)^{1/2}(1+\left\Vert X^1\right\Vert _{p}).  \label{G'}
\end{equation}
\item[$\mathbf{S}(\beta)$] [\textbf{Stochastic Sewing} with $\beta>0$] There exists $C_{sew}>0$  such that  \red{for every $s\leq r\leq t$}%
\begin{align}
& \left\Vert \widehat{\mathcal{E}}_{s,r,t}^{2}(X)-\widehat{\mathcal{E}}%
_{s,r,t}^{1}(X)\right\Vert _{p}\leq C_{sew}(t-s)\big(\left\Vert
X^{1}-X^{2}\right\Vert _{p}+(1+\Vert X\Vert _{p})(t-s)^{\beta }\big)
\label{V7} \\
& \left\Vert \widetilde{\mathcal{E}}_{s,r,t}^{2}(X)-\widetilde{\mathcal{E}}%
_{s,r,t}^{1}(X)\right\Vert _{p}\leq C_{sew}(t-s)^{1/2}\big(\left\Vert
X^{1}-X^{2}\right\Vert _{p}+(1+\Vert X\Vert _{p})(t-s)^{\beta }\big),
\label{V8}
\end{align}
\red{where} $\|X\|_p=\|X^1\|_p+\|X^2\|_p$ for $X=(X^1,X^2)$. 
\end{itemize}


\begin{definition}\label{def_betacoupling}
  A family of operators $\mathcal{X}_{s,r,t}(X)$, $s\leq r\leq t$ which satisfies
  the properties $\mathbf{R}$, $\mathbf{G}$ and $\mathbf{S}(\beta)$ will be called a $\beta$-coupling family of operators for $\Theta$. 
  
  If such a representation holds, we will say that the family $\Theta_{s,t}:\mathcal{P}_p(H)\to\mathcal{P}_p(H)$ satisfies the $\beta$-sewing property.
  \end{definition}
  
  For $X_{0}\in L_{t_{0}}^{p}(H)\times L_{t_{0}}^{p}(H)$, a partition $\pi =\{t_{0}<t_{1}<....<t_{N}\}\subset [0,T]$ and some $r_k\in [t_k,t_{k+1})$, $0\le k\le N-1$,  
  we define the random variables $X_{k}$, $k=1,\ldots ,N$, recursively as 
  \begin{equation}\label{def_Xk}
  X_{k}=\mathcal{X}_{t_{k-1},r_{k-1},t_{k}}(X_{k-1}).
  \end{equation}%
  
  We observe that from representation property {\bf R}, we get from (\ref{def_Thetapi}) that $X^1_N\sim \Theta^{\pi}_{t_0,t_N}(\mathcal{L}(X_0^1))$ and $X^2_N\sim \Theta^{\pi'}_{t_0,t_N}(\mathcal{L}(X_0^2))$, where $\pi'=\{t_0\le r_0<t_1\le r_1<...<r_{N-1}<t_N\}$.  
  The scheme $(X_k,0\le k \le N)$ will allow us to analyse the distance between these two distributions. Before this, we start by giving  an upper bound of the marginal laws. For this, it is enough to focus on the first marginal. This is why the growth assumption {\bf G} is only written on the first marginal.   

\begin{lemma}
\label{BASIC1} Suppose that $\mathbf{R}$ and $\mathbf{G}$ hold. Let $T>0$ and $t_N\le T$.
Then, there exists a constant $A_T \in \R_+$  depending on $T$, $p$ and $L_{sew}$ such that  
\begin{eqnarray}
1+\|\Theta^\pi_{t_0,t_N}(\mu)\|_p &\le& A_T (1+\|\mu\|_p).\label{majo_Lambda}
\end{eqnarray}%
\end{lemma}
The proof of this lemma is postponed to~Appendix~\ref{AppendixlemBASIC1}.

\begin{definition}\label{def_subpartition}
  Let  $\pi
  =\{t_{0}<t_{1}<....<t_{N}\}$ be a partition of $[t_0,t_N]$.   
  \red{A subpartition of $\pi$ is a partition $\pi^{\prime }=\{t^{\prime
  }_0<\dots<t^{\prime }_{N^{\prime }}\}$ with $N'\ge N$ such that $t'_0=t_0$, $t'_{N'}=t_N$ and 
  there exist $0<i_1<\dots<i_{N-1}<N'$ such that $t'_{i_k}=t_{k}$
  for $k\in \{1,\dots, N-1\}$.}
  
  A subpartition of $\pi$ is a simple subpartition of $\pi$ if it is such that 
  $\pi ^{\prime }=\{t_{0}\leq r_{0}<t_{1}<....<t_{N-1}\leq r_{N-1}<t_{N}\}$
  and has at most $2N+1$ elements ($r_i$ and $t_i$ are ``merged'' if $r_i=t_i$%
  ).
  \end{definition}

The main result in this section is the following stochastic sewing lemma, whose proof is
postponed in Appendix~\ref{AppendixlemBASIC2}.  
\begin{lemma}
\label{sewing} (Wasserstein sewing lemma) Assume that $\Theta_{s,t}:\mathcal{P}_p(H) \to \mathcal{P}_p(H)$ is a family satisfying the $\beta$-sewing property for some $\beta>0$ (i.e.
there exists a joint representation of $\Theta_{s,t}$ and $\Theta_{r,t}\circ \Theta_{s,r}$ for which {\bf R}, {\bf G} and {\bf S}($\beta$) hold). Let \red{$\pi =\{t_{0}<t_{1}<....<t_{N}\}$ be a partition of $[0,T]$, $T>0$,}   and $\pi'$ a simple subpartition of $\pi$. 
Then, for any
$\mu,\nu \in \mathcal{P}_p(H)$, we have 
\begin{equation}
  W_{p}(\Theta _{t_{0},t_{N}}^{\pi }(\mu ),\Theta _{t_{0},t_{N}}^{\pi ^{\prime
  }}(\nu ))\leq C_T W_{p}(\mu ,\nu
  )+C_T(1+\|\mu\|_p+\|\nu\|_p)(t_{N}-t_{0})^{1/2}\left\vert \pi \right\vert
  ^{\beta },  \label{V9'}
  \end{equation}%
  where $C_T \in \mathbb{R}_+$ is a constant depending on $T$, $p$, $L_{sew}$ and $C_{sew}$.
  In particular, taking $\pi ^{\prime }=\pi ,$ we get the asymptotic Lipschitz
  property%
  \begin{equation}
  W_{p}(\Theta _{t_{0},t_{N}}^{\pi }(\mu ),\Theta _{t_{0},t_{N}}^{\pi }(\nu
  ))\leq C_T W_{p}(\mu ,\nu )+C_T (1+\|\mu\|_p+\|\nu\|_p)(t_{N}-t_{0})^{1/2}\left\vert
  \pi \right\vert ^{\beta }.  \label{V9"}
  \end{equation}%
  
\end{lemma}
\begin{remark}
  The ``sewing lemma" denomination was introduced in the seminal papers by Feyel, de La Pradelle \cite
  {[FP]} and Feyel, de La Pradelle, Mokobodzki \cite{[FPM]}. Here, the property $\mathbf{S}(\beta)$ is the appropriate version of the
  stochastic sewing property introduced recently by Le in\cite[Theorem 2.1, Eq.~(2.8) and (2.9)]{[Le]}, which allows to take advantage of the martingale property.
 Let us note that for $\beta>1/2$, $\nu=\mu \in \mathcal{P}_p(H)$, $\pi=\{s<t\}$ and $\pi'=\{s<r<t\}$, we get from~\eqref{V9'} the local estimate
 $$W_p(\Theta_{s,t}(\mu),\Theta_{r,t}\circ \Theta_{s,r}(\mu))\le C_T(1+2\|\mu\|_p)(t-s)^{1/2+\beta},$$
which is the ``classical" sewing property. It allows to construct a unique limit flow, see~\cite[Lemma 2.1 and Proposition 2.4]{[AB]}. 
The interest of the stochastic sewing lemma is thus to handle cases with $\beta \in (0,1/2]$.  
In contrast with Le~\cite{[Le]}, where the stochastic sewing assumption is directly stated at the level of random variables,  we work here on probability measures and  need to represent them by random variables, which is made with the operator $\mathcal{X}_{s,r,t}$ and the sequence~\eqref{def_Xk}. The assumptions made on $\mathcal{X}_{s,r,t}$ are local in time, but cannot be translated directly in terms of $\Theta_{s,t}$, contrary to the classical sewing lemma. The estimate~\eqref{V9'} is not local in time  but brings on partitions, which allows to take advantage of the martingale argument that underlies the stochastic sewing lemma. This estimate is then sufficient to construct then a limit flow.
  \end{remark}

\subsection{Construction of the flow $\theta$}

We consider a family of operators $\Theta_{s,t}:\mathcal{P}_p(H)\to \mathcal{P}_p(H)$, $0\le s\le t$, that satisfies the $\beta$-sewing property (see Definition~\ref{def_betacoupling}). We aim at
 constructing a flow $\theta _{s,t}$ obtained as the limit of the schemes built with $\Theta_{s,t}$. 

To make this construction, we will use dyadic partitions of step $\frac{1}{2^{n}}$. We first introduce some notation. \red{For $n \in \N$, we denote $\DD_n=\{k2^{-n},k \in \N\}$. We also define $\DD=\cup_{n\in \N} \DD_n$ the set of nonnegative dyadic numbers.
For $t\geq 0$ and $n\in \N$ we then define $\eta_{\DD_n}(t)= \frac{\lfloor 2^n t\rfloor}{2^n}$ so that
\begin{eqnarray*}
&& \eta _{\DD_n}(t)=\frac{k}{2^{n}}\quad \text{if}\quad \frac{k}{2^{n}}%
\leq t<\frac{k+1}{2^{n}}. \\
\end{eqnarray*}%
} \red{For $s<t$, we consider the following partition of step $\frac{1}{2^{n}}$: } 
\begin{equation}
\pi _{n}=\pi _{n}(s,t)=\{\eta_{\DD_n}(s)< \eta_{\DD_n}(s)+2^{-n}< ...<\eta_{\DD_n}(t)-2^{-n}<\eta_{\DD_n}(t)\}.  \label{grid}
\end{equation}%
Notice that this partition does not cover exactly the interval $[s,t]$  but $[\eta
_{\DD_n}(s),\eta _{\DD_n}(t)]$.
Then, we denote 
\begin{equation*}
\Theta _{s,t}^{n}=\Theta_{\eta_{\DD_n}(s),\eta_{\DD_n}(t)}^{\pi _{n}(\red{s,t})}
=\Theta
_{\eta_{\DD_n}(t)-2^{-n},\eta_{\DD_n}(t)}\circ \dots \circ \Theta
_{\eta_{\DD_n}(s),\eta_{\DD_n}(s)+2^{-n}}.
\end{equation*}%
It is a scheme with a constant time step equal to $2^{-n}$ between $\eta_{\DD_n}(s)$ and $\eta_{\DD_n}(t)$.
As a consequence of Lemma~\ref{sewing}, we obtain the next result.

\begin{proposition}
Suppose that the family $\Theta_{s,t}$ satisfies the $\beta$-sewing property. Let $T>0$. If $s,t\in  \DD_{n}$  with $s<t\le T$, then 
\begin{equation}
W_{p}(\Theta _{s,t}^{n}(\mu ),\Theta _{s,t}^{n+1}(\nu ))\leq C_TW_{p}(\mu
,\nu )+C_T (1+\|\mu\|_p+\|\nu\|_p)(t-s)^{1/2}2^{-\beta n} . \label{V9}
\end{equation}
where $C_T\in \R_+$  is a  constant depending on $T$, $p$, $L_{sew}$ and $C_{sew}$. Moreover, 
\begin{equation}
W_{p}(\Theta _{s,t}^{n}(\mu ),\Theta _{s,t}(\mu ))\leq C'_T (1+\|\mu\|_p)(t-s)^{1/2},  \label{V10bis}
\end{equation}%
where $C'_T\in \R_+$  is a  constant depending on $T$, $p$, $L_{sew}$, $C_{sew}$ and $\beta$.
\end{proposition}

\textbf{Proof. } Since $s,t \in \DD_{n}$, $\pi_{n+1}(s,t)$ is a simple subpartition
of $\pi_n(s,t)$. We can thus apply Lemma \ref{sewing}: the inequality (\ref{V9}) is a direct consequence of (\ref{V9'}).

We now check (\ref{V10bis}). For $t-s>1$, we  get $W_{p}(\Theta _{s,t}^{n}(\mu ),\Theta _{s,t}(\mu )) \leq  \| \Theta _{s,t}^{n}(\mu ) \|_p+ \| \Theta _{s,t}(\mu )\|_p \le 2 A_T \| \mu \|_p$ by Lemma~\ref{BASIC1}. 

We now assume $t-s\le 1$ and introduce some
notation: for a grid $\pi =\{t_{0}<...<t_{k}<...<t_{N}\}$ we denote $\lambda
(\pi )=\pi \setminus \{t_{2k+1},1\leq 2k+1<N\}$. So in order to get the
partition $\lambda(\pi)$ we drop out the points of the partition $\pi$ that have an odd index. Then, we iteratively define $\lambda ^{m}(\pi )=\lambda(\lambda ^{m-1}(\pi ))$%
. We use this for the grid $\pi _{n}(s,t)$ defined in (\ref{grid}). We
notice that if $s,t\in \DD_{n} $ then, for $m\leq n,$ we have the following bound on the time-step
$$\left\vert
\lambda ^{m}(\pi _{n}(s,t))\right\vert \leq \frac{1}{2^{n-m}}.$$ We now
observe that 
\red{$\#\lambda(\pi)=\lfloor\frac{\# \pi}2\rfloor +1$} is a
\red{nondecreasing} function of $\#\pi$. We also observe that if $\#\pi=2^{n+1}+1$
for some $n\in \mathbb{N}$, then $\#\lambda(\pi)=2^{n}+1$. Therefore, using that $%
t\leq s+1$ \red{and $\#\pi_n(s,s+1)=2^{n}+1$}, we get 
\begin{equation*}
\# \lambda^n(\pi_n(s,t))\le \# \lambda^n(\pi_n(s,s+1))=2.
\end{equation*}
This gives that $\lambda ^{n}(\pi
_{n}(s,t))=\{s,t\}$ so that $\Theta _{s,t}^{\lambda ^{n}(\pi _{n}(s,t))}(\mu
)=\Theta _{s,t}(\mu ).$ Consequently, we get by using the triangle inequality and then (\ref{V9'}), 
\begin{eqnarray*}
W_{p}(\Theta _{s,t}^{n}(\mu ),\Theta _{s,t}(\mu )) &=&W_{p}(\Theta
_{s,t}^{\pi _{n}(s,t)}(\mu ),\Theta _{s,t}^{\lambda ^{n}(\pi _{n}(s,t))}(\mu
)) \\
&\leq &\sum_{i=1}^{n}W_{p}(\Theta _{s,t}^{\lambda ^{i-1}(\pi _{n}(s,t))}(\mu
),\Theta _{s,t}^{\lambda ^{i}(\pi _{n}(s,t))}(\mu )) \\
&\leq &C_T (1+2\|\mu\|_p)(t-s)^{1/2}\sum_{i=1}^{n}\left\vert \lambda ^{i-1}(\pi
_{n}(s,t))\right\vert ^{\beta } \\
&\leq &C_T  (1+2\|\mu\|_p)(t-s)^{1/2}\sum_{i=1}^{n}\frac{1}{2^{(n+1-i)%
\beta }}\\
&\leq& \frac{C_T}{2^\beta -1}  (1+2\|\mu\|_p)(t-s)^{1/2}.
\square\end{eqnarray*}%

In order to state the next result, we need to introduce the following uniform continuity property on the diagonal:
\begin{itemize}
\item[$\mathbf{C}$] \textbf{[Continuity]} There exists an increasing function $%
\varphi :[0,\infty)\rightarrow [0,\infty)$ such that $\lim_{\delta
\rightarrow 0}\varphi (\delta )=0$ and, for every $s<t$ and every $\mu \in 
\mathcal{P}_{p}(H)$%
\begin{equation}  \label{cont}
W_{p}(\Theta _{s,t}(\mu ),\mu )\leq (1+\|\mu\|_p)\varphi (t-s).
\end{equation}
\end{itemize}

\begin{theorem}
\label{existence}Suppose that the family $\Theta_{s,t}:\mathcal{P}_p(H) \to \mathcal{P}_p(H)$, $0\le s<t$, satisfies the $\beta$-sewing property (Def.~\ref{def_betacoupling}) and $\mathbf{C}$.
There exists $\theta _{s,t}:\mathcal{P}_{p}(H)\rightarrow \mathcal{P}_{p}(H)$, $0\le s\le t$,
which is a flow (i.e. $\theta _{s,t}=\theta _{r,t}\circ\theta _{s,r}$ for $s\leq r \leq t$) and such that $W_p(\Theta^n_{s,t}(\mu),\theta_{s,t}(\mu))\to_{n\to \infty} 0$ for any $s<t$, $\mu \in \mathcal{P}_p(H)$. It satisfies for any $T>0$:
\begin{itemize}
\item the linear growth property
\begin{equation}\label{lin_growth}
\|\theta_{s,t}(\mu)\|_p\le A_T(1+\|\mu\|_p),  \ 0 \le s <t \le T,
\end{equation}   
with $A_T\in \R_+$  depending on $T$, $p$ and $L_{sew}$,
\item  the following Lipschitz property 
\begin{equation}
W_{p}(\theta _{s,t}(\mu ),\theta _{s,t}(\nu ))\leq C_T W_{p}(\mu ,\nu ), \ 0 \le s <t \le T
\label{V13}
\end{equation}%
where $C_T\in \R_+$  is depending on $T$, $p$, $L_{sew}$ and $C_{sew}$,
\item $(s,t)\mapsto \theta_{s,t}(\mu)$ is continuous uniformly in $(\mathcal{P}_p(H),W_p)$
w.r.t. $\mu$ on $\{0\le s <t \le T\}$.
\end{itemize}
\end{theorem}

\textbf{Proof. Step 1 [construction of $\theta_{s,t}$ for dyadic times $s<t$
]} Let $s<t$ be fixed such that $s,t\in \DD_{n}=\{\frac{k}{%
2^{n}},k\in \N\}$ for some $n\in \mathbb{N}$. Then,  for $i\ge n$, $\pi_{i+1}(s,t)$ is a simple subpartition of $\pi_{i}(s,t)$. \red{Taking} $m_2\ge m_1\ge n$,  we get by using (\ref{V9}) 
\begin{eqnarray}
W_{p}(\Theta _{s,t}^{m_1}(\mu ),\Theta _{s,t}^{m_2}(\mu )) &\leq
&\sum_{i=m_1}^{m_2-1}W_{p}(\Theta _{s,t}^{i}(\mu ),\Theta _{s,t}^{i+1}(\mu )) \notag \\
&\leq &C_T  (1+2\|\mu\|_p)(t-s)^{1/2}\sum_{i=m_1}^{m_2-1}2^{-\beta i} \notag\\
&\le& \frac{2^\beta C_T}{2^\beta -1}  (1+2\|\mu\|_p)(t-s)^{1/2}2^{-m_1\beta }. \label{ineg_step1}
\end{eqnarray}%
Therefore,  the sequence $\Theta_{s,t}^{m}(\mu
)$, $m\ge n$ is Cauchy and we define $\theta _{s,t}(\mu )=\lim_{m \rightarrow
\infty }\Theta _{s,t}^{m}(\mu )$ (note that $(\mathcal{P}_p(H),W_p)$ is
complete by \cite[Theorem 6.18]{Villani}). Since $n$ is arbitrary, we have defined 
$\theta _{s,t}(\mu )$ for every dyadic number $s,t\in \DD:=\cup
_{n=1}^{\infty }\DD_{n}$ with $s<t$.

We now show~(\ref{lin_growth}) and~(\ref{V13}) when $s,t \in \DD$, with $s<t$.  Let $n$ be such that $s,t\in \DD_{n}$. From~\eqref{majo_Lambda}, we easily get~(\ref{lin_growth}). From the asymptotic Lipschitz property~(\ref{V9"}), we get that $W_{p}(\theta_{s,t}(\mu ),\theta_{s,t}(\nu )) =\lim_{m \to \infty}  W_{p}(\Theta _{s,t}^{m}(\mu ),\Theta _{s,t}^{m}(\nu )) \le C_T W_p(\mu,\nu)$.

\smallskip \textbf{Step 2 [flow property on dyadic times]} Let us prove the
flow property $\theta _{s,t}=\theta _{r,t} \circ \theta _{s,r}$ for $s,r,t\in \DD$ such that $s<r<t$. Let $s,r,t\in \DD_{n}$ and $m\geq n$. One has $\Theta _{s,t}^{m}\red{(\mu)}=\Theta
_{r,t}^{m} \circ \Theta _{s,r}^{m}\red{(\mu)}.$ We write the triangle inequality
\begin{eqnarray*}
&&W_{p}(\theta _{r,t} \circ \theta _{s,r}(\mu ),\Theta _{r,t}^{m} \circ \Theta
_{s,r}^{m}(\mu )) \\
&\leq &W_{p}(\theta _{r,t}(\theta _{s,r}(\mu )),\Theta _{r,t}^{m}(\theta
_{s,r}(\mu )))+W_{p}(\Theta _{r,t}^{m}(\theta _{s,r}(\mu )),\Theta
_{r,t}^{m}(\Theta _{s,r}^{m}(\mu ))).
\end{eqnarray*}%
By the very construction of $\theta$ the first term vanishes as $m\to \infty$. From the
asymptotic Lipschitz continuity property (\ref{V9"}), the second term is upper bounded by 
\begin{equation*}
C_T W_{p}(\theta _{s,r}(\mu ),\Theta _{s,r}^{m}(\mu ))+C_T
(1+\|\theta _{s,r}(\mu )\|_p+\|\Theta _{s,r}^{m}(\mu )\|_p)(r-s)^{1/2}
2^{-m\beta }.
\end{equation*}%
By Lemma~\ref{BASIC1} $\|\theta _{s,r}(\mu )\|_p+\|\Theta _{s,r}^{m}(\mu )\|_p\le 2A_T (1+\|\mu\|_p)$  is bounded in $m$, so the second term goes to $0$  as $m\rightarrow \infty$. Thus, we get $W_{p}(\theta _{r,t} \circ \theta _{s,r}(\mu ),\Theta _{r,t}^{m} \circ \Theta
_{s,r}^{m}(\mu )) \to_{m\to \infty} 0$.

Finally, using the triangular inequality and $%
\Theta_{s,t}^{n}(\mu )=\Theta_{r,t}^{n} \circ \Theta _{s,r}^{n}(\mu )$, we get 
\begin{eqnarray*}
&&W_{p}(\theta _{s,t}(\mu ),\theta _{r,t}\circ \theta _{s,r}(\mu )) \\&\leq
&W_{p}(\Theta _{s,t}^{n}(\mu ),\theta _{s,t}(\mu )) +W_{p}(\theta
_{r,t}\circ \theta _{s,r}(\mu ),\Theta _{r,t}^{n}\circ\Theta _{s,r}^{n}(\mu ))
\rightarrow_{n\to \infty} 0.
\end{eqnarray*}

\textbf{Step 3 [continuity on dyadic times and extension]} Using (\ref%
{V10bis}) and (\ref{cont}), for $t-s\leq 1,s,t\in \DD$, we get
\begin{eqnarray*}
  W_{p}(\Theta _{s,t}^{n}(\mu ),\mu) &\leq &W_{p}(\Theta _{s,t}^{n}(\mu ),\Theta _{s,t}(\mu
))+W_{p}(\Theta _{s,t}(\mu ),\mu ) \\
&\leq & {C_T^{\prime
}(1+\|\mu\|_p)}\sqrt{t-s}+(1+\|\mu\|_p)\varphi (t-s).
\end{eqnarray*}%
Letting $n\to \infty$, we obtain%
\begin{equation}\label{cont_flow}
W_{p}(\theta _{s,t}(\mu ),\mu )\leq {C_T^{\prime
}(1+\|\mu\|_p)}\sqrt{t-s}+(1+\|\mu\|_p)\varphi (t-s).
\end{equation}%
Recall that we have the flow property $\theta _{s,t+h}=\theta
_{t,t+h} \circ \theta_{s,t}$ for $s,t,t+h \in \DD$ such that $s<t<t+h$. Applying the previous inequality, we get%
\begin{eqnarray*}
W_{p}(\theta_{s,t+h}(\mu ),\theta _{s,t}(\mu )) &=&W_{p}(\theta_{t,t+h}(\theta _{s,t}(\mu )),\theta _{s,t}(\mu )) \\
&\leq &{C_T^{\prime
}(1+\|\theta _{s,t}(\mu )\|_p)}\sqrt{h}+(1+\|\theta _{s,t}(\mu )\|_p)\varphi (h)\\
&\leq &{A_T C_{T}^{\prime}}(1+\|\mu\|_p)(\sqrt{h}+\varphi(h))\rightarrow_{h\to 0} 0.
\end{eqnarray*}%
  Similarly, we get for $s,s+h^{\prime },t\in \DD$ such that 
$s<s+h^{\prime }\le t$ 
\begin{align*}
W_{p}(\theta _{s,t}(\mu ),\theta _{s+h^{\prime },t}(\mu ))&=W_{p}(\theta
_{s+h^{\prime },t}(\theta_{s,s+h^{\prime }}(\mu) ),\theta _{s+h^{\prime
},t}(\mu ))\le C_T W_{p}(\theta_{s,s+h^{\prime }}(\mu),\mu ) \\
&\le C_T C^{\prime }_T(1+\|\mu\|_p)(\sqrt{h^{\prime }}+
\varphi(h^{\prime })) \rightarrow_{h^{\prime }\to 0} 0,
\end{align*}
by using~\eqref{V13} (that has been proved in Step 1 for dyadic numbers) and then~\eqref{cont_flow}.

These estimates allow us to extend by continuity $(s,t)\red{\mapsto} \theta
_{s,t}(\mu )$ from $(s,t)\in \DD\times \DD$ to any $s<t$. Namely, we have for $m\le n$, 
\begin{align*}
  W_p(\theta_{\eta_{\DD_m}(s),\eta_{\DD_m}(t)}(\mu),\theta_{\eta_{\DD_n}(s),\eta_{\DD_n}(t)}(\nu))\le  (A_T+C_T) C^{\prime }_T (1+\|\mu\|_p) (2^{-m/2}+\varphi(2^{-m})),
\end{align*}
and $\theta_{s,t}(\mu)$ is the limit of the Cauchy sequence $\theta_{\eta_{\DD_n}(s),\eta_{\DD_n}(t)}(\mu)$. Note that this is also the limit of $\Theta^n_{s,t}(\mu)$, since $$W_p(\Theta^n_{s,t}(\mu),\theta^n_{\eta_{\DD_n}(s),\eta_{\DD_n}(t)}(\mu))\le \frac{2^\beta C_T}{2^\beta -1} (1+2\|\mu\|_p) T^{1/2}2^{-n \beta }$$ by~(\ref{ineg_step1}).

Clearly, $\theta$ thus constructed has the flow property, and the properties~(\ref{lin_growth}) and~(\ref{V13}) hold. We also  have 
\begin{align*}W_p(\theta_{s,t+h}(\mu),\theta_{s+h',t}(\mu))&\le W_p(\theta_{s,t+h}(\mu),\theta_{s,t}(\mu))+W_p(\theta_{s,t}(\mu),\theta_{s+h',t}(\mu))\\
  &\le   (A_T+C_T) C^{\prime }_T (1+\|\mu\|_p)\left(\sqrt{h}+
  \varphi(h)+\sqrt{h^{\prime }}+
  \varphi(h^{\prime }) \right),
\end{align*}
which gives the uniform continuity. $\square $

\subsection{Speed of convergence}

 We have defined in Theorem~\ref{existence} the flow~$\theta$ as the limit of the scheme on the dyadic grid and then, we have extended it by continuity on~$\mathbb{R}_+$. Here, we give a convergence result for any partition, as well as the speed of convergence. This result has an interest in itself, and will enables us to state a uniqueness result.

\begin{theorem}
\label{main} Let $\Theta_{s,t}:\mathcal{P}_p(H)\to \mathcal{P}_p(H)$ be a family of operators satisfying the $\beta$-sewing property (Def.~\ref{def_betacoupling}) and $\mathbf{C}$. Let $\theta$ be the flow given by Theorem~\ref{existence}. Let $T>0$. Then, for any $0\le s<t \le T$, there exists a constant $C'_T\in {%
\mathbb{R}}_+$ that only depends on  $T$, $p$, $L_{sew}$, $C_{sew}$ and $\beta$ such that for
every partition $\pi $ of $[s,t]$%
\begin{equation}
W_{p}(\Theta _{s,t}^{\pi }(\mu ),\theta_{s,t} (\mu ))\leq C'_T (1+\|\mu\|_p)%
(t-s)^{1/2}\left\vert \pi \right\vert ^{\beta}.  \label{m1}
\end{equation}%
In particular, for the ``trivial partition" we have $t-s=\left\vert \pi
\right\vert $ so we get%
\begin{equation}
W_{p}(\Theta _{s,t}(\mu ),\theta_{s,t} (\mu ))\leq C'_T (1+\| \mu\|_p)%
(t-s)^{(1/2+\beta )}.  \label{m3}
\end{equation}
\end{theorem}

\textbf{Proof. Step 1 [dyadic partitions]} In a first step we fix $n \ge 0$
and we assume that 
\begin{equation*}
\pi =\{t_{0}<t_{1}<...<t_{m}\}\quad \text{ with }\quad t_{i}=\frac{k_{i}}{2^{n}}, \  t_m\le T.
\end{equation*}
Our goal is to prove the following estimate 
\begin{equation}  \label{estim_dyadic}
W_{p}(\Theta _{t_{0},t_{m}}^{\pi }(\mu ),\Theta _{t_{0},t_{m}}^{\pi _{n}(t_0,t_m)}(\mu ))\leq C'_T (1+\|\mu\|_p) (t_{m}-t_{0})^{1/2}\left\vert \pi \right\vert ^{\beta },
\end{equation}
where $\pi_n(t_0,t_n)$ is defined by~\eqref{grid}.
 The specificity of the partition $%
\pi $ is that the times $t_{i}$ are in the dyadic grid of step $\frac{1}{%
2^{n}}$: in a second step we will get rid of this restriction.

Now, we add some intermediary points between $t_{i}=\frac{k_{i}}{%
2^{n}}$ and $t_{i+1}=\frac{k_{i+1}}{2^{n}}$ as follows. We denote%
\begin{equation*}
d(k_{i},k_{i+1})=\left\lfloor \frac{k_{i}+k_{i+1}}{2}\right\rfloor
\end{equation*}%
and we define 
\begin{equation*}
\gamma (\pi )=\{\frac{k_{0}}{2^{n}}\leq \frac{d(k_{0},k_{1})}{2^{n}}<\frac{%
k_{1}}{2^{n}}\leq ...\leq \frac{k_{m-1}}{2^{n}}\leq \frac{d(k_{m-1},k_{m})}{%
2^{n}}<\frac{k_{m}}{2^{n}}\}
\end{equation*}

\begin{remark}
Roughly speaking, one obtains the partition $\gamma (\pi )$ by adding in
between each $\frac{k_{i}}{2^{n}}$ and $\frac{k_{i+1}}{2^{n}}$ the ``middle
point" $\frac{d(k_{i},k_{i+1})}{2^{n}}.$ But, if there is no precise ``middle
point" belonging to the dyadic grid $\DD_{n},$ we take the point to the left
(given by the integer part). Moreover, notice that we may have $%
k_{i+1}=k_{i}+1$ and in this case $d(k_{i},k_{i+1})=k_{i}.$ This means that
we do not add a supplementary point if there is no free place in between $%
\frac{k_{i}}{2^{n}}$ and $\frac{k_{i+1}}{2^{n}}.$
\end{remark}

Note that $\gamma(\pi)$ is a simple subpartition of $\pi$ with elements in
the dyadic grid $\DD_{n}$. We have $\pi=\gamma(\pi)$ if, and only if $\pi=\pi
_{n}(t_{0},t_{m})$ defined in (\ref{grid}), which is the finest partition of $%
[t_0,t_m]$ on the dyadic grid. Otherwise, we have $\#\gamma(\pi)>\#\pi$. The
partition $\pi _{n}(t_{0},t_{m})$ is thus a subpartition of $\pi$ and of $%
\gamma^i(\pi)$ for $i\ge 1$. Therefore, the sequence $(\# \gamma^i(\pi))_{i\ge
1}$ is nondecreasing and bounded $\# \gamma^i(\pi)\le k_m-k_0+1$, and it is
increasing until it reaches this bound. We set 
\begin{align*}
q& =\inf \{i\ge 0, \gamma^i(\pi)= k_m-k_0+1 \}=\inf \{i\ge 0, \gamma^i(\pi)=
\gamma^{i+1}(\pi) \} \\
&=\inf \{i\ge 0, \gamma^i(\pi)= \pi_n(t_0,t_m) \}.
\end{align*}

We now prove~\eqref{estim_dyadic}. The case $q=0$ is trivial ($%
\pi=\pi_n(t_0,t_m)$) and we assume $q\ge 1$. For $1\le i\le q$, $%
\gamma^{i}(\pi)$ is a simple subpartition of $\gamma^{i-1}(\pi)$. We can
thus use~(\ref{V9'}) and get 
\begin{eqnarray}
W_{p}(\Theta _{t_{0},t_{m}}^{\pi }(\mu ),\Theta _{t_{0},t_{m}}^{\gamma
^{q}(\pi )}(\mu )) &\leq & \sum_{i=1}^{q} W_{p}(\Theta
_{t_{0},t_{m}}^{\gamma ^{i-1}(\pi )}(\mu ) ,\Theta _{t_{0},t_{m}}^{\gamma
^{i}(\pi )}(\mu ))  \notag \\
&\leq & C_T (1+2\|\mu\|_p)(t_{m}-t_{0})^{1/2}\sum_{i=1}^{q}\left\vert \gamma ^{{i-1}}(\pi )\right\vert ^{\beta }.
\label{estim_dyadic_2}
\end{eqnarray}
Now, we observe that $|\gamma^i(\pi)|=\nu_i2^{-n}$ for some natural number $%
\nu_i\ge 1$ ($\nu_0=\max_{1\le \ell \le m}k_{\ell}-k_{\ell-1}$), and we have
the induction formula: $\nu_i=\frac{\nu_{i-1}}{2}$ if $\nu_{i-1}$ is even
and $\nu_i=\frac{\nu_{i-1}+1}{2}$ if $\nu_{i-1}$ is odd. In all cases, we
have $\nu_i\le \frac{\nu_{i-1}+1}{2}$ and thus 
\begin{equation*}
\nu_i- 1 \le 2^{-i} (\nu_0-1).
\end{equation*}
Therefore, we have $\sum_{i=1}^{q}\left\vert \gamma ^{{i-1}}(\pi
)\right\vert ^{\beta }\le 2^{-n\beta} \sum_{i=0}^{q-1}
(1+2^{-i}(\nu_0-1))^\beta$. Now, by definition of $q$, we have $%
\nu_{q-1}\ge 2$ and thus $q-1\le \frac{\log(\nu_0-1)}{\log(2)}$. We get 
\begin{align*}
\frac 1 {\nu_0^\beta} \sum_{i=0}^{q-1} (1+2^{-i}(\nu_0-1))^\beta
&\le 2^{(\beta-1)^+}\sum_{i=0}^{ \frac{\log(\nu_0-1)}{\log(2)}} \left(
\frac 1 {\nu_0^\beta} + 2^{-i\beta}\right) \\
&\le 2^{(\beta-1)^+}\left( \sup_{\nu \ge 2 }\frac{1+\frac{\log(\nu-1)}{%
\log(2)}}{ \nu^\beta }+\frac{1}{1-2^{-\beta}}\right):=\bar
C_\beta
\end{align*}
Therefore we get $%
\sum_{i=1}^{q}\left\vert \gamma ^{{i-1}}(\pi )\right\vert ^{\beta }\le
\bar C_\beta 2^{-n\beta}\nu_0^\beta=
\bar C_\beta|\pi|^\beta$. Thus~\eqref{estim_dyadic} holds with $C'_T=2\bar C_\beta C_T$  by using~%
\eqref{estim_dyadic_2}.

\smallskip

\textbf{Step 2 [general partitions]} Take a general partition 
\begin{equation*}
\pi =\{t_{0}<...<t_{m}\}.
\end{equation*}%
We recall that $\eta _{\DD_n}(t)=\frac{k}{2^{n}}$ if $\frac{k}{2^{n}}\leq t<%
\frac{k+1}{2^{n}}$. We take $n$ sufficiently large (i.e. $2^{-n}<\min_{1\le
\ell\le m}t_\ell-t_{\ell-1}$) so that $\eta _{\DD_n}(t_{0})<...<\eta _{\DD_n}(t_{m})$%
, and we define the partition 
\begin{equation*}
\pi_{n}=\{\eta _{\DD_n}(t_{0})<...<\eta _{\DD_n}(t_{m})\}.
\end{equation*}%
We have by the triangular inequality 
\begin{eqnarray*}
&&W_{p}(\Theta _{t_{0},t_{m}}^{\pi }(\mu ),\theta _{t_{0},t_{m}}(\mu )) \\
&\leq &W_{p}(\Theta _{t_{0},t_{m}}^{\pi }(\mu ),\Theta _{\eta
_{\DD_n}(t_{0}),\eta _{\DD_n}(t_{m})}^{\pi _{n}}(\mu ))+W_{p}(\Theta _{\eta
_{\DD_n}(t_{0}),\eta _{\DD_n}(t_{m})}^{\pi _{n}}(\mu ),\Theta _{\eta
_{\DD_n}(t_{0}),\eta _{\DD_n}(t_{m})}^{\pi _{n}(t_{0},t_{m})}(\mu )) \\
&&+W_{p}(\Theta _{\eta _{\DD_n}(t_{0}),\eta _{\DD_n}(t_{m})}^{\pi
_{n}(t_{0},t_{m})}(\mu ),\theta _{\eta _{\DD_n}(t_{0}),\eta _{\DD_n}(t_{m})}(\mu
))+W_{p}(\theta _{\eta _{\DD_n}(t_{0}),\eta
_{\DD_n}(t_{m})}(\mu ),\theta _{t_{0},t_{m}}(\mu )) \\
&=:&\sum_{i=1}^{4}a_{i}(n).
\end{eqnarray*}%
From~(\ref{estim_dyadic}), we get 
\begin{equation*}
a_{2}(n)=W_{p}(\Theta_{\eta _{\DD_n}(t_{0}),\eta _{\DD_n}(t_{m})}^{\pi _{n}}(\mu
),\Theta _{\eta _{\DD_n}(t_{0}),\eta _{\DD_n}(t_{m})}^{\pi _{n}(t_{0},t_{m})}(\mu
))\leq  C'_T (1+\|\mu\|_p)(\eta_n(t_m)-\eta_n(t_0))^{1/2}%
\left\vert \pi _{n}\right\vert ^{\beta }
\end{equation*}%
and since $|\pi_n|\leq |\pi|+2^{-n}$, we have 
\begin{equation*}
\limsup_{n\to\infty}a_{2}(n)\leq  C'_T (1+\|\mu\|_p)(t_m-t_0)^{1/2}\left\vert \pi \right\vert^{\beta }.
\end{equation*}

By Lemma \ref{lemma-app1}, we get
\begin{equation*}
a_1(n)\leq C_T (1+\|\mu\|_p) \Big(\varphi\Big(\frac 1{2^n}\Big)+\Big(%
t_m-t_0+\frac 1{2^n}\Big)^{1/2}\Big(|\pi|+\frac 1{2^n}\Big)^{\beta}%
\Big),
\end{equation*}
which gives 
\begin{equation*}
\limsup_{n\to\infty}a_1(n)\leq C_T(1+\|\mu\|_p)(t_m-t_0)^{1/2}|\pi|^{%
\beta}.
\end{equation*}
Moreover we have $\lim_{n}a_{3}(n)=0$ from~(\ref{ineg_step1}) and by the continuity of $\theta$ obtained in Theorem~\ref{existence} gives $\lim_{n}a_{4}(n)=0$. We therefore get  $$W_{p}(\Theta _{t_{0},t_{m}}^{\pi }(\mu ),\theta _{t_{0},t_{m}}(\mu )) \le (C'_T+C_T) (1+\|\mu\|_p)
(t_m-t_0)^{1/2}\left\vert \pi \right\vert^{\beta }$$ so the proof is completed. $%
\square $

\section{Application to stochastic jump equations}\label{Sec_jump}

We start by presenting the stochastic equation that we want to study in this section.  We consider some Hilbert separable space $H$ (in particular, the standard example is $H=\R^{d})$, and our aim is to solve the
following equation on $H:$%
\begin{eqnarray}
X_{s,t} &=&X+\int_{s}^{t}\int_{H\times E}c(v,z,X_{s,r-},\mathcal{L}%
(X_{s,r})))\widetilde{N}_{\mathcal{L}(X_{s,r})}(dv,dz,dr)  \label{eq1} \\
&&+\int_{0}^{t}\int_{H\times E}\sigma (v,z,X_{s,r},\mathcal{L}(X_{s,r}))W_{%
\mathcal{L}(X_{s,r})}(dv,dz,dr)  \notag \\
&&+\int_{s}^{t}b(X_{s,r},\mathcal{L}(X_{s,r}))ds, \ t\in [s,+\infty),  \notag
\end{eqnarray}
with
\begin{equation*}
c,\sigma :H\times E\times H\times \mathcal{P}_{2}(H)\rightarrow H,\quad
b:H\times \mathcal{P}_{2}(H)\rightarrow H.
\end{equation*}
Let us precise what are the objects in the above equation. This is an equation with values in $H$: so $X_{s,t}\in H$ and $\mathcal{L}(X_{s,r})$ is the distribution of $X_{s,r}$, which is assumed to be probability measure in $\mathcal{P}_{2}(H)$.

Here $(E,\mathcal{E})$ is a Lusin space and $\nu $ is a $\sigma$-finite
measure on $E$. Moreover, for a family of probability measures $\eta _{s,t}\in \mathcal{P}_{2}(H)$, $s<t$,  such that $(s,t)\mapsto \eta_{s,t}$ is measurable, and for each fixed $s$
 we define $N_{\mathcal{\eta}}$ to be a Poisson point measure on $%
H\times E\times (s,\infty )$ with compensator 
\begin{equation*}
\widehat{N}_{\mathcal{\eta }}(dv,dz,dr)=\eta_{s,r}(dv)\nu (dz)1_{(s,\infty
)}(r)dr.
\end{equation*}%
We define the compensated measure  $\widetilde{N}_{\mathcal{\eta }}=N_{\mathcal{\eta }}-\widehat{N%
}_{\mathcal{\eta }}$ \red{(we use the notation from Ikeda and Watanabe \cite{[IW]})}. We denote 
\begin{equation*}
\mathcal{F}_{s,t}^{N}=\ \sigma (N_{\mathcal{\eta }}(A):A\in \mathcal{B}%
(H)\times \mathcal{E}\times \mathcal{B}(s,t))
\end{equation*}%
The important fact concerning the stochastic integral with respect to $%
\widetilde{N}_{\mathcal{\eta }}(dv,dz,dr)$ is that it is a linear operator
such that for a process $\varphi :(s,\infty )\times H\times E\rightarrow H$
which is previsible with respect to $(\mathcal{F}_{s,t}^{N})_{t\geq s}$, one
has the isometry property%
\begin{eqnarray*}
\E\left(\left\vert \int_{s}^{t}\int_{H\times E}\varphi (r,v,z)\widetilde{N}%
_{\mathcal{\eta }}(dv,dz,dr)\right\vert _{H}^{2}\right) &=&\E\int_{s}^{t}\int_{H\times E}\left\vert \varphi (r,v,z)\right\vert
_{H}^{2}\mathcal{\eta }_{s,r}(dv)\nu (dz)dr.
\end{eqnarray*}%
The stochastic integral makes sense only if the last term is finite.

Moreover, $W_{\mathcal{\eta }}(dv,dz,dr)$ is a martingale measure  on $H\times
E\times (s,\infty )$ of intensity 
\red{$$
\widehat{W}_{\mathcal{\eta }}(dv,dz,dr)=\mathcal{\eta }_{s,r}(dv)\nu (dz)1_{(s,+\infty)}(r)dr
$$ }
(see e.g. El Karoui and \red{M\'el\'eard}~\cite[Theorem I-4]{ElKM}).
We denote by 
\begin{equation*}
\mathcal{F}_{s,t}^{W}=\ \sigma (W_{\mathcal{\eta }}(A):A\in \mathcal{B}%
(H)\times \mathcal{E}\times \mathcal{B}(s,t)),
\end{equation*}%
its natural filtration and we set $\mathcal{G}_{s,t}=\mathcal{F}_{s,t}^{W}\vee \mathcal{F}_{s,t}^{N},t\geq
s$ completed with sets of probability zero.

Again, the stochastic integral with respect to $%
W_{\mathcal{\eta }}(dv,dz,dr)$ is a linear operator such that for
a process $\varphi :(s,\infty )\times H\times E\rightarrow H$ which is
previsible with respect to $(\mathcal{F}_{s,t}^{W})_{t\geq s}$ one has the
isometry property 
\begin{eqnarray*}
&&\E\left(\left\vert \int_{s}^{t}\int_{\red{H}\times E}\varphi (r,v,z)W_{\mathcal{\eta }%
}(dv,dz,dr)\right\vert _{H}^{2}\right) \\
&=&\E\int_{s}^{t}\int_{\red{H}\times E}\left\vert \varphi (r,v,z)\right\vert
_{H}^{2}\mathcal{\eta }_{s,r}(dv)\nu (dz)dr.
\end{eqnarray*}%
The stochastic integral makes sense if the last term is finite. 

We are now able to give a sense to a solution of Equation~(\ref{eq1}): this is the existence of a process~$(X_{s,t})_{t\ge 0}$, a Poisson random measure~$N_\eta$ and a martingale measure $W_\eta$ with $\eta_{s,t}=\mathcal{L}(X_{s,t})$ such that~(\ref{eq1}) hold for any $t\ge s$.

Before going further, we introduce \red{the} ``canonical" Poisson process and martingale measure that we will use to construct our coupling operator. Namely, let $N^{can}$ (resp $W^{can}$) be a Poisson point measure (resp. martingale measure) on $(0,1)\times E \times (0,\infty)$ with intensity  $\mathbf{1}_{(0,1)}(w)dw \nu(dz) \mathbf{1}_{(0,+\infty)}(r)dr$. We assume that $N^{can}$ and $W^{can}$ are independent, and consider $U$ an independent uniform random variable on $(0,1)$.  Let $\mathcal{F}^{can}_t$ be the natural filtration associated to $U$ and the random measures $N^{can}$ and $W^{can}$, completed with sets of probability zero. The filtered probability space \red{satisfies} the atomless assumption {\bf A}, and we will use intensively Lemma~\ref{lem_dist} to represent the Euler scheme for~\eqref{eq1}  with $N^{can}$ and $W^{can}$ and to construct a suitable coupling operator.

\subsection{Assumptions and main result}

Before stating our main result of the section, we present our assumptions on the coefficients
\begin{equation*}
c,\sigma :H\times E\times H\times \mathcal{P}_{2}(H)\rightarrow H,\quad
b:H\times \mathcal{P}_{2}(H)\rightarrow H.
\end{equation*}%

We assume that they verify the following Lipschitz continuity property.
\begin{itemize}
\item[\textbf{L}] \textbf{[Lipschitz property]} For every $v,v^{\prime }\in
H,x,x^{\prime }\in H$ and every $\rho ,\rho ^{\prime }\in \mathcal{P}_{2}(H)$
\begin{eqnarray}
&&\left\vert c(v,z,x,\rho )-c(v^{\prime },z,x^{\prime },\rho ^{\prime
})\right\vert _{H}+\left\vert \sigma (v,z,x,\rho )-\sigma (v^{\prime
},z,x^{\prime },\rho ^{\prime })\right\vert _{H}  \label{a1} \\
&\leq &C(z)(\left\vert v-v^{\prime }\right\vert _{H}+\left\vert x-x^{\prime
}\right\vert _{H}+W_{2}(\rho ,\rho ^{\prime }))  \notag
\end{eqnarray}%
with 
\begin{equation}
\int_{E}C^{2}(z)\nu (dz)=:\overline{C}<\infty  \label{a1'}
\end{equation}%
and 
\begin{equation}
\left\vert b(x,\rho )-b(x^{\prime },\rho ^{\prime })\right\vert _{H}\leq
C(\left\vert x-x^{\prime }\right\vert _{H}+W_{2}(\rho ,\rho ^{\prime })).
\label{a2}
\end{equation}
\end{itemize}

We assume in addition that 
\begin{equation}  \label{a2'}
\exists \rho_0 \in \mathcal{P}_2, \ \int_E |c(0,z,0,\rho_0)|_H^2+|\sigma(0,z,0,\rho_0)|_H^2
\nu(dz)<\infty.
\end{equation}
Then, we may assume without loss of generality that $$|b(0,\rho_0)|\le C \text{ and } 
|c(0,z,0,\rho_0)|+|\sigma(0,z,0,\rho_0)|\le C(z)$$ (take $\tilde{C}(z)=\max(C(z),|c(0,z,0,\rho_0)|_H+|\sigma(0,z,0,\rho_0)|_H)$
instead of $C(z)$ otherwise), so that we get from~(\ref{a1}) and (\ref{a2})
the sublinear properties: 
\begin{align}
& \left\vert c(v,z,x,\rho ) \right\vert _{H}+\left\vert \sigma (v,z,x,\rho )
\right\vert _{H} \leq C(z)(1+ \left\vert v\right\vert _{H}+\left\vert x
\right\vert _{H}+W_{2}(\rho ,\rho_0)) \label{sublin_csig}\\
& \left\vert b(x,\rho ) \right\vert _{H}\leq C(1+\left\vert x\right\vert
_{H}+W_{2}(\rho ,\rho_0)).\label{sublin_b}
\end{align}

\begin{theorem}
  \label{jump} Suppose that (\ref{a1}), (\ref{a1'}), (\ref{a2}) and (%
  \ref{a2'}) hold. We have:
  \begin{itemize}
    \item Existence. For any $X\sim\mu\in \mathcal{P}_2(H)$, there exists a solution to Equation (\ref{eq1}), i.e. there exists a process $(X_{s,t})_{t\ge s}$, a Poisson random measure $N_{\mathcal{L}(X_{s,t})}$ and a martingale measure $W_{\mathcal{L}(X_{s,t})}$ such that~(\ref{eq1}) holds. 
    \item Uniqueness of marginal laws. $\theta_{s,t}(\mu):= \mathcal{L}(X_{s,t})$ is continuous with respect to~$t$ for~$W_2$, and any solution $\overline{X}_{s,t}$  such that $t\mapsto  \mathcal{L}(\overline{X}_{s,t})$ is continuous for~$W_2$ satisfies $\mathcal{L}(\overline{X}_{s,t})=\theta_{s,t}(\mu)$.
    \item Properties. The map $\theta_{s,t}:\mathcal{P}_2(H)\to\mathcal{P}_2(H)$ enjoys the flow property $\theta _{r,t}\circ \theta _{s,r}(\mu )=\theta _{s,t}(\mu )$ for $s<r<t$, and we have the following Lipschitz property 
  \begin{equation}\label{Lip_jump}
  \forall T>0,\exists C_T, \forall 0<s<t<T, \ W_{2}(\theta _{s,t}(\mu ),\theta _{s,t}(\nu ))\leq C_T W_{2}(\mu ,\nu  ).
  \end{equation}%
  Moreover, \red{there exists a family of operators  $\Theta_{s,t}:\mathcal{P}_2(H)\to \mathcal{P}_2(H)$, $t\geq s$, such that the following holds:} there exists $C'_T\in \R_+$ such that for every partition $\pi$ of $(s,t)$ with $t\le T$, we have%
  \begin{equation}\label{CV_jump}
  W_{2}(\mathcal{L}(X_{s,t}(X)),\Theta _{s,t}^{\pi }(\mathcal{L}(X)))\leq C'_T (1+\|X\|_2)
  \sqrt{t-s} \left\vert \pi \right\vert^{1/2}.
  \end{equation}
\end{itemize}
  
  \end{theorem}

To prove Theorem~\ref{jump}, we introduce the Euler scheme \red{$\Theta_{s,t}(\mu)$, $\mu\in \mathcal{P}_2(H)$, appearing in \eqref{CV_jump}}, that we will use in order to construct and approximate the solution of this equation. Given $\mu \in \mathcal{P}_2(H)$ we
consider a Poisson point measure $N_{\mu }$ of intensity $\mu (dv)\times \nu
(dz)\times dr$ and a martingale measure $W_{\mu }$ with the same intensity (this corresponds to the \red{processes already} defined with $\eta_{s,t}=\mu$ being constant).
Then, we take a random variable $X\sim \mu $ which is independent of $N_{\mu} $ and of $W_{\mu }$ and we define%
\begin{eqnarray}
X^{Eul}_{s,t}(X) &=&X+\int_{s}^{t}\int_{H\times E}c(v,z,X,\mathcal{\mu })%
\widetilde{N}_{\mathcal{\mu }}(dv,dz,\red{dr})  \label{euler} \\
&&+\int_{s}^{t}\int_{H\times E}\sigma (v,z,X,\mathcal{\mu })W_{\mathcal{\mu 
}}(dv,dz,\red{dr})  \notag \\
&&+\int_{s}^{t}b(X,\mathcal{\mu })\red{dr},  \notag
\end{eqnarray}%
and
\begin{equation}
\Theta _{s,t}(\mu ):=\mathcal{L}(X^{Eul}_{s,t}(X)).  \label{Euler}
\end{equation}%
We think to $\Theta _{s,t}(\mu )$ as to an one step Euler scheme starting
from $\mu$. Moreover, for a partition $\pi =\{t_{0}<t_{1}<....<t_{N}\}$ we recall~(\ref{def_Thetapi})
\begin{equation*}
\Theta _{t_0,t_N}^{\pi }(\mu )=\Theta _{t_{N-1},t_{N}}\circ ...\circ
\Theta _{t_{0},t_{1}}(\mu ),
\end{equation*}%
which is the law of the Euler scheme on the time grid~$\pi$ starting at time $t_0$ from the distribution $\mu$.

\smallskip
\textbf{Sketch of the proof.} The proof of Theorem~\ref{jump} is split into three steps\red{, developed in Sections \ref{Step1}, \ref{Step2} and \ref{Step3}}.

\medskip

\textbf{Step 1.} We prove that $\Theta$ satisfies the $\beta$-sewing property and 
$\mathbf{C}$. Then we use Theorem~\ref{existence} in order to construct the flow $\theta_{s,t}$ associated to $\Theta _{s,t}.$

\medskip

\textbf{Step 2.} Let us consider two families of maps $\mathcal{\rho }_{s,t},\mathcal{\rho}'_{s,t}:\mathcal{P}_2(H)\to \mathcal{P}_2(H)$, $0\le s<t$ 
such that for any $\mu$, $(s,t)\mapsto \red{\mathcal{\rho
}_{s,t}}(\mu)$ and $(s,t)\mapsto \red{\mathcal{\rho
}'_{s,t}}(\mu)$ are  measurable.  
For a given $\mu \in \mathcal{P}_2(H)$, we consider $X\sim \mu$ and the compensated random Poisson and Wiener measures $\widetilde{N}_{\rho (\mu)}$ and $%
W_{\rho(\mu)}$ that are independent, and we consider the (standard) Stochastic Differential Equation with jumps: 
\begin{eqnarray}
X_{s,t}^{\rho } &=&X+\int_{s}^{t}\int_{H\times E}c(v,z,X_{s,r-}^{\rho },%
\mathcal{\rho }_{s,r}(\mu ))\widetilde{N}_{\rho(\mu)}(dv,dz,dr)
\label{eq1bis} \\
&&+\int_{s}^{t}\int_{H\times E}\sigma (v,z,X_{s,r-}^{\rho },\rho
_{s,r}(\mu ))W_{\rho(\mu) }(dv,dz,dr)  \notag \\
&&+\int_{s}^{t}b(X_{s,r},\mathcal{\rho }_{s,r}(\mu ))ds.  \notag
\end{eqnarray}
We define similarly $X_{s,t}^{\rho'}$ and prove a stability property which, roughly speaking, \red{gives}
a bound of $W_{2}(\mathcal{L}(X_{s,t}^{\rho }(X)),\mathcal{L}%
(X_{s,t}^{\rho ^{\prime }}(X)))$ in terms of $W_{2}(\rho _{s,t},\rho
_{s,t}^{\prime })$.

\medskip

\textbf{Step 3.} We prove existence of a solution to Equation~(\ref{eq1}) in the following way: we consider the SDE~(\ref{eq1bis}) with $\rho_{s,t}=\theta _{s,t}$ and $\mu =\mathcal{L}(X),$ where $\theta$ is the flow obtained in Step 1. For such an equation, the existence of solutions is standard. The delicate point  is to prove that the law of
the solution of this equation is indeed given by $\theta _{s,t}(\mu )$. This is
done by using the Euler scheme approximation and the stability property
proved at Step 2. Uniqueness follows also by the stability property from Step 2.

\subsection{Step 1: construction of $\theta _{s,t}$}\label{Step1}

In view of applying Theorem~\ref{main}, we first construct the coupling operator.  Namely, we consider $\Theta _{s,t}$ defined in (\ref{Euler}) and we give a suitable 
stochastic representation (coupling operator) of $\Theta _{s,t}$ and $\Theta
_{r,t}\circ \Theta _{s,r}$. We fix $s<r<t$ and $X=(X^{1},X^{2})\in
L_{s}^{2}(H)\times L_{s}^{2}(H)$ and we denote $\mathcal{L}(X^{i})=\mu^{i}.$ Our goal is to construct a random variable $\mathcal{X}_{s,r,t}(X^{1},X^{2})=(\mathcal{X}^1_{s,r,t}(X^{1},X^{2}),\mathcal{X}^2_{s,r,t}(X^{1},X^{2}))$ such that $$\mathcal{L}(\mathcal{X}^1_{s,r,t}(X^{1},X^{2}))=\Theta_{s,t}(\mu^1) \text{ and }   \mathcal{L}(\mathcal{X}^2_{s,r,t}(X^{1},X^{2}))=\Theta_{r,t} \circ \Theta_{s,r}(\mu^2),$$ 
so that {\bf R} holds. We also want that {\bf S}($\beta$) and {\bf C} hold for this representation. 

Let $\Pi _{s}(dv_{1},dv_{2})$ be the probability measure on $H\times H$ which is a $W_{2}$-optimal coupling of $\mu^{1}$ and $\mu^{2}$ and we take $\tau_{s}:(0,1)\rightarrow H\times H$ such that 
\begin{equation*}
\int_{\red{H\times H}}\varphi (v_{1},v_{2}) \Pi
_{s}(dv_{1},dv_{2})=\int_{0}^{1}\varphi (\tau _{s}(w))dw.
\end{equation*}
Such a coupling exists by~\cite[Theorem 4.1]{Villani}. Similarly, we take $\Pi _{r}(dv_{1},dv_{2})$ to be the
probability measure on $H\times H$ which is an optimal coupling of $\mu^{1}$ and $\Theta_{s,r}(\mu^{2})$, and
we take $\tau _{r}:(0,1)\rightarrow H\times H$ such that 
\begin{equation*}
\int_{H\times H}\varphi (v_{1},v_{2})\Pi
_{r}(dv_{1},dv_{2})=\int_{0}^{1}\varphi (\tau _{r}(w))dw.
\end{equation*}%
Notice that, by construction, 
\begin{eqnarray*}
\int_{0}^{1}|\tau _{s}^{1}(w)-\tau _{s}^{2}(w)|_{H}^{2}dw &=&W_{2}^{2}(\mu
^{1},\mu^{2}) \\
\int_{0}^{1}|\tau _{r}^{1}(w)-\tau _{r}^{2}(w)|_{H}^{2}dw &=&W_{2}^{2}(\mu
^{1},\Theta _{s,r}(\mu^{2})) .\\
\end{eqnarray*}

Now we are able to give the $L^{2}$ representation. Let us recall that $N^{can}$ (resp. $W^{can}$) is a Poisson point
process (resp. martingale measure) \red{on $(0,1)\times E\times (0,+\infty)$} of intensity $dw\times \nu
(dz)\times dr$.
Then, we construct%
\begin{eqnarray*}
\mathcal{X}_{s,r}^{2}(X^{1},X^{2}) &=&X^{2}+b(X^{2},\mu ^{2})(r-s) \\
&&+\int_{s}^{r}\int_{(0,1)\times E}c(\tau _{s}^{2}(w),z,X^{2},\mu^{2}) 
\widetilde{N}^{can}(dw,dz,du) \\
&&+\int_{s}^{r}\int_{(0,1)\times E}\sigma (\tau_{s}^{2}(w),z,X^{2},\mu^{2})W^{can}(dw,dz,du)
\end{eqnarray*}%
and%
\begin{eqnarray*}
\mathcal{X}_{s,r,t}^{2}(X^{1},X^{2}) &=&\mathcal{X}_{s,r}^{2}(X^{1},X^{2})+b(%
\mathcal{X}_{s,r}^{2}(X^{1},X^{2}),\Theta _{s,r}(\mu^{2}))(t-r) \\
&&+\int_{r}^{t}\int_{(0,1)\times E}c(\tau _{r}^{2}(w),z,\mathcal{X}%
_{s,r}^{2}(X^{1},X^{2}),\Theta _{s,r}(\mu^{2})) \widetilde{N}^{can}(dw,dz,du) \\
&&+\int_{r}^{t}\int_{(0,1)\times E}\sigma (\tau _{r}^{2}(w),z,\mathcal{X}%
_{s,r}^{2}(X^{1},X^{2}),\Theta _{s,r}(\mu^{2})) W^{can}(dw,dz,du)
\end{eqnarray*}%
It is easy to check (see Appendix~\ref{append_schemes}) that the law of $\mathcal{X}_{s,r}^{2}(X^{1},X^{2})$ is $
 \Theta _{s,r}(\mu^{2})$. 
In the same way, the law of $\mathcal{X}%
_{s,r,t}^{2}(X^{1},X^{2})$ is $\Theta _{r,t}\circ \Theta _{s,r}(\mu
_{s}^{2}).$ Next we construct
\begin{eqnarray*}
\mathcal{X}_{s,r}^{1}(X^{1},X^{2}) &=&X^{1}+b(X^{1},\mu^{1})(r-s) \\
&&+\int_{s}^{r}\int_{(0,1)\times E}c(\tau _{s}^{1}(w),z,X^{1},\mu^{1}) 
\widetilde{N}^{can}(dw,dz,du) \\
&&+\int_{s}^{r}\int_{(0,1)\times E}\sigma (\tau _{s}^{1}(w),z,X^{1},\mu^{1}) W^{can}(dw,dz,du)
\end{eqnarray*}%
and%
\begin{eqnarray*}
\mathcal{X}_{s,r,t}^{1}(X^{1},X^{2}) &=&\mathcal{X}%
_{s,r}^{1}(X^{1},X^{2})+b(X^{1},\mu^{1})(t-r) \\
&&+\int_{r}^{t}\int_{(0,1)\times E}c(\tau_{r}^{1}(w),z,X^{1},\mu^{1}) 
\widetilde{N}^{can}(dw,dz,du) \\
&&+\int_{r}^{t}\int_{(0,1)\times E}\sigma (\tau _{r}^{1}(w),z,X^{1},\mu^{1}) W^{can}(dw,dz,du)
\end{eqnarray*}%
So in both intervals $(s,r)$ and $(r,t)$ we keep $X^{1}$ and $\mu^{1}$
but we change $\tau_{s}^{1}$ with $\tau _{r}^{1}$ (which appears from the
optimal coupling of $\mu^{1}$ with $\Theta_{s,r}(\mu^{2})$). We check again easily that the law of $\mathcal{X}_{s,r,t}^{1}(X^{1},X^{2})$ is $\Theta _{s,t}(\mu
^{1}).$

\begin{remark}
Notice that the definition of $\Theta _{s,t}(\mu )$ is given by means of the probabilistic representation (\ref{euler}). However, the stochastic representation given by $\mathcal{X}_{s,r,t}(X^{1},X^{2})$ is different from a simple composition of $\Theta$. The main point is that we change the coupling between the Poisson and martingale measures that drive the coarse and the fine scheme, by changing from $\tau_s$ to $\tau_r$. This is crucial to get then the $\beta$-sewing property. 
\end{remark}

\begin{proposition}\label{prop_RSC}
  Under the hypotheses of Theorem~\ref{jump},  $\Theta$ defined by \eqref{Euler} satisfies \red{the $\beta$-sewing property with $\beta=1/2$}
%
%
  and  $\mathbf{C}$. 
\end{proposition}
{\bf Proof.} To prove the $1/2$-sewing property, we show that the stochastic representation $\mathcal{X}_{s,r,t}(X^{1},X^{2})$ satisfies \textbf{R}, $\mathbf{G}$ and \textbf{S}(1/2).

The representation property \textbf{R} has been checked while $\mathbf{G}$ is an easy consequence of~(\ref{sublin_csig}) and~(\ref{sublin_b}). We show~\textbf{S}(1/2). We have from~(\ref{def_E_hat_tilde}) 
\begin{equation*}
\widetilde{\mathcal{E}}_{s,r,t}^{2}(X^{1},X^{2})-\widetilde{\mathcal{E}}%
_{s,r,t}^{1}(X^{1},X^{2})=I_{s,r}^{N}+I_{s,r}^{W}+I_{r,t}^{N}+I_{r,t}^{W}+J_{r,t}
\end{equation*}%
with%
\begin{eqnarray*}
I_{s,r}^{N} &=&\int_{s}^{r}\int_{(0,1)\times E}(c(\tau
_{s}^{2}(w),z,X^{2},\mu^{2})-c(\tau _{s}^{1}(w),z,X^{1},\mu^{1})) 
\widetilde{N}^{can}(dw,dz,du) \\
I_{s,r}^{W} &=&\int_{s}^{r}\int_{(0,1)\times E}(\sigma (\tau
_{s}^{2}(w),z,X^{2},\mu^{2})-\sigma (\tau _{s}^{1}(w),z,X^{1},\mu^{1})) W^{can}(dw,dz,du),
\end{eqnarray*}%
and 
\begin{eqnarray*}
I_{r,t}^{N} &=&\int_{r}^{t}\int_{(0,1)\times E}(c(\tau _{r}^{2}(w),z,%
\mathcal{X}_{s,r}^{2}(X^{1},X^{2}),\Theta_{s,r}(\mu^{2}))-c(\tau
_{r}^{1}(w),z,X^{1},\mu^{1})) \widetilde{N}^{can}(dw,dz,du) \\
I_{r,t}^{W} &=&\int_{r}^{t}\int_{(0,1)\times E}(\sigma (\tau _{r}^{2}(w),z,%
\mathcal{X}_{s,r}^{2}(X^{1},X^{2}),\Theta_{s,r}(\mu^{2}))-\sigma (\tau
_{r}^{1}(w),z,X^{1},\mu^{1})) W^{can}(dw,dz,du)
\end{eqnarray*}%
and finally

\begin{equation*}
J_{r,t}=\bigg(b(\mathcal{X}_{s,r}^{2}(X^{1},X^{2}),\Theta_{s,r}(\mu^{2}))-\E\big(b(\mathcal{X}_{s,r}^{2}(X^{1},X^{2}),\Theta_{s,r}(\mu^{2}))\big| \mathcal{F%
}_{s} \big)\bigg)(t-r)
\end{equation*}%
Using the isometry property first and the Lipschitz property \textbf{L} next, we get%
\begin{eqnarray*}
\Vert I_{s,r}^{N}\Vert _{2}^{2} &=&\int_{s}^{r}\int_{(0,1)\times E}\E \big(%
|c(\tau _{s}^{2}(w),z,X^{2},\mu^{2})-c(\tau _{s}^{1}(w),z,X^{1},\mu^{1})|_{H}^{2}\big)dw \nu (dz)du \\
&\leq &3 \overline{C}(r-s)\left(\int_{0}^{1}\left\vert \tau _{s}^{2}(w)-\tau
_{s}^{1}(w)\right\vert _{H}^{2}dw+\E\left\vert X^{2}-X^{1}\right\vert
_{H}^{2}+W_{2}^{2}(\mu _{s}^{2},\mu _{s}^{1})\right).
\end{eqnarray*}%
Notice that $W_{2}^{2}(\mu _{s}^{2},\mu _{s}^{1})\leq \E \left\vert
X^{1}-X^{2}\right\vert _{H}^{2}$ and 
\begin{equation*}
\int_{0}^{1}\left\vert \tau _{s}^{2}(w)-\tau _{s}^{1}(w)\right\vert
_{H}^{2}dw=W_{2}^{2}(\mu _{s}^{1},\mu _{s}^{2})\leq \E\left\vert
X^{1}-X^{2}\right\vert _{H}^{2}
\end{equation*}%
so that%
\begin{equation*}
\Vert I_{s,r}^{N}\Vert _{2}^{2}\leq 9 \overline{C}(r-s)\Vert X^{2}-X^{1}\Vert_{2}^{2}.
\end{equation*}%
The same calculus gives 
\begin{equation*}
\Vert I_{s,r}^{W}\Vert _{2}^{2}\leq 9\overline{C}(r-s)\E\Vert
X^{2}-X^{1}\Vert _{2}^{2}.
\end{equation*}%
In the same way, we get
\begin{eqnarray*}
\left\Vert I_{r,t}^{N}\right\Vert _{2}^{2} &=&\int_{r}^{t}\int_{(0,1)\times
E}\E\big(|c(\tau _{r}^{2}(w),z,\mathcal{X}_{s,r}^{2}(X^{1},X^{2}),\Theta_{s,r}(\mu^{2}))-c(\tau _{r}^{1}(w),z,X^{1},\mu^{1})|_{H}^{2}\big)dw \nu (dz)
du \\
&\leq &3 \overline{C}(t-r) \bigg(\int_{0}^{1}\left\vert \tau _{r}^{2}(w)-\tau
_{r}^{1}(w)\right\vert _{H}^{2}dw+ \E\left\vert \mathcal{X}%
_{s,r}^{2}(X^{1},X^{2})-X^{1}\right\vert _{H}^{2} \\
&&+W_{2}^{2}(\Theta_{s,r}(\mu^{2}),\mu^{1})\bigg) \\
&=&3\overline{C}(t-r)\bigg(\E\left\vert \mathcal{X}_{s,r}^{2}(X^{1},X^{2})-X^{1}%
\right\vert _{H}^{2}+2W_{2}^{2}(\Theta_{s,r}(\mu^{2}),\mu _{r}^{1})\bigg)\\
&\le &9\overline{C}(t-r) \E\left\vert \mathcal{X}_{s,r}^{2}(X^{1},X^{2})-X^{1}%
\right\vert _{H}^{2}
\end{eqnarray*}%
since $\mathcal{X}_{s,r}^{2}(X^{1},X^{2})\sim \Theta _{s,r}(\mu
_{s}^{2}) = \mu^2_r$.
This gives%
\begin{eqnarray*}
\Vert I_{r,t}^{N} \Vert _{2}^{2} &\leq &18\overline{C}(t-r)( \Vert X^{2}-X^{1} \Vert_{2}^{2} + \Vert \mathcal{X}_{s,r}^{2}(X^{1},X^{2})-X^{2}\Vert_2^2)
\end{eqnarray*} 
By using the isometry and the sublinear properties~(\ref{sublin_csig}--\ref{sublin_b}) on $b$, $c$, $\sigma$, we
get 
\begin{align}
&\Vert \mathcal{X}_{s,r}^{2}(X^{1},X^{2})-X^{2} \Vert_2^2
\notag \\
&\leq 3 \E[|b(X^2,\mu^2)|_H^2](r-s)^2 +3 \E\left[\int_s^r\int_{(0,1) \times
E }|c(\tau^2_s(w),z,X^2,\mu^2)|_H^2dw \nu (dz)dt \right]  \notag \\
& \ +3 \E\left[\int_s^r\int_{(0,1) \times E
}|\sigma(\tau^2_s(w),z,X^2,\mu^2)|_H^2dw \nu (dz)dt \right]  \notag \\
&\leq \tilde{C} (1+\|X^2\|_2)^2 ((r-s)+(r-s)^2) \leq \tilde{C}(1+T)(1+\|X^2\|_2)^2
 (r-s).  \label{bound_increment}
\end{align}
We finally get the existence of a constant $C_T>0$ such that 
\begin{equation*}
\Vert I_{r,t}^{N} \Vert _{2}^2\le C_T (t-r) \left( \|X^2-X^1\|_2^2+
(1+\|X^2\|_2)^2 (r-s) \right).
\end{equation*}

The same estimate holds for $I_{r,t}^{W}$. Last, we  write
$$J_{r,t}=(t-r)\left(b(\mathcal{X}_{s,r}^{2}(X^{1},X^{2}),\Theta_{s,r}(\mu^{2}))-b(X^{2},\mu^{2}) - \E[b(\mathcal{X}_{s,r}^{2}(X^{1},X^{2}),\Theta _{s,r}(\mu^{2}))-b(X^{2},\mu^{2}) |\mathcal{F}_s]\right)$$
and thus we get by \textbf{L}, 
\begin{eqnarray*}
\E(|J_{r,t}|_{H}^{2}) &\leq &2 \E\big(|b(\mathcal{X}_{s,r}^{2}(X^{1},X^{2}),\Theta _{s,r}(\mu^{2}))-b(X^{2},\mu^{2})|_{H}^{2}\big)(t-r)^{2} \\
&\leq &C \Vert\mathcal{X}_{s,r}^{2}(X^{1},X^{2})-X^{2}\Vert_{2}^{2}%
(t-r)^{2},
\end{eqnarray*}%
that can be upper bounded by~(\ref{bound_increment}). So we have proved that 
\begin{equation*}
\Vert \widetilde{\mathcal{E}}_{s,r,t}^{2}(X^{1},X^{2})-\widetilde{\mathcal{E}}_{s,r,t}^{1}(X^{1},X^{2})\Vert _{2}^{2}\leq C_T
(t-s)\left( \Vert X^{2}-X^{1}\Vert_{2}^{2} +  (1+\|X^2\|_2)^2 (t-s) \right).
\end{equation*}

On the other hand, we have from~(\ref{def_E_hat_tilde}) 
\begin{align*}
&\widehat{\mathcal{E}}_{s,r,t}^{2}(X^{1},X^{2})-\widehat{\mathcal{E}}%
_{s,r,t}^{1}(X^{1},X^{2})= (b(X^{2},\mu^{2})-b(X^{1},\mu
^{1}))(r-s) \\
& +\left( \E\left((b(\mathcal{X}_{s,r}^{2}(X^{1},X^{2}),\Theta _{s,r}(\mu^{2})))-b(X^{2},\mu^{2})+b(X^{2},\mu^{2})-b(X^{1},\mu^{1})\big|{\mathcal{F}_{s}}\right)\right) (t-r).
\end{align*}%
Using the Lipschitz property \textbf{L} and~(\ref{bound_increment}), we get 
\begin{equation*}
\Vert\widehat{\mathcal{E}}_{s,r,t}^{2}(X^{1},X^{2})-\widehat{\mathcal{E}}%
_{s,r,t}^{1}(X^{1},X^{2})\Vert_{2}^{2}\leq C_{T}(t-s)^{2}\left( \Vert
X^{1}-X^{2}\Vert_{2}^{2}+(1+\|X^2\|_2)^2(t-s)\right)
\end{equation*}%
which means that the property $\mathbf{S}(\beta)$  holds with $\beta =1/2.$

We finally prove the continuity property $\mathbf{C}$. We come back to (\ref{euler}) and get with the same calculation as in~(\ref{bound_increment})
\begin{equation*}
W_{2}^{2}(\Theta _{s,t}(\mu ),\mu )\leq \E\left\vert X_{s,t}(X)-X\right\vert
_{H}^{2}\leq C(1+\|\mu\|_2)^2\left((t-s)^2+(t-s)\right)
\end{equation*}%
so $\mathbf{C}$\ holds with $\varphi (u)=\sqrt{C}(u+\sqrt{u})$, $u\ge 0$.$\square$


\medskip

As a consequence of Proposition~\ref{prop_RSC}, we can apply Theorem~\ref{existence}, and get the existence of a flow $\theta_{s,t}$ that is continuous for $W_2$, satisfies~\eqref{Lip_jump} and also~\eqref{CV_jump} by Theorem~\ref{main}. To get Theorem~\ref{CV_jump}, it remains to prove the existence and the uniqueness results, for which we need a stability result.

\subsection{Step 2: a stability property.}\label{Step2}

In this section we consider a family of applications $\rho _{s,t}:\mathcal{P}%
_{2}(H)\rightarrow \mathcal{P}_{2}(H)$ such that for each $s\ge 0$, $t\mapsto \rho_{s,t}(\mu)$ is continuous. Let  $c,\sigma $ and 
$b$ be coefficients satisfying (\ref{a1}), (\ref{a1'}), (\ref{a2}) and (\ref%
{a2'}). For $s\ge 0$ fixed, we consider the stochastic differential equation, for $t\ge s$, 
\begin{eqnarray}
X_{s,t}^{\rho } &=&X+\int_{s}^{t}\int_{H\times E}c(v,z,X_{s,r-}^{\rho },%
\mathcal{\rho }_{s,r}(\mu ))\widetilde{N}_{\mathcal{\rho }}(dv,dz,dr)
\label{stab1} \\
&&+\int_{s}^{t}\int_{H\times E}\sigma (v,z,X_{s,r}^{\rho },\rho _{s,r}(\mu
)) W_{\mathcal{\rho }}(dv,dz,{dr})  \notag \\
&&+\int_{s}^{t}b(X_{s,r}^{\rho },\mathcal{\rho }_{s,r}(\mu )){dr}.  \notag
\end{eqnarray}
Here, $X\sim \mu \in \mathcal{P}_2(H)$, and the compensated Poisson measure $\widetilde{N}_{\mathcal{\rho }}$ and the martingale measure 
$W_{\rho}$ have  the intensity ${\rho_{s,r}(\mu)(dv)\nu(dz)dr}$. Following~\cite[Theorem 9.1 p.~245]{[IW]},  there exists a unique c\`adl\`ag solution to this equation.

Then, we fix $n\in \N$ and we consider the time grid $t_{k}^{n}=s+\frac{k}{n}$. For $t_{k}^{n}\leq r<t_{k+1}^{n}$, we denote $\eta_{n}(r)=t_{k}^{n}.$ Then, we associate the discretized equation 
\begin{eqnarray}
X_{s,t}^{\rho ,n} &=&X+\int_{s}^{t}\int_{H\times E}c(v,z,X_{s,\red{\eta _{n}(r)-}}^{\rho ,n},\mathcal{\rho }_{s,\eta _{n}(r)}(\mu ))\widetilde{N}_{%
\mathcal{\rho }}(dv,dz,dr)  \label{stab2} \\
&&+\int_{s}^{t}\int_{H\times E}\sigma (v,z,X_{s,\red{\eta _{n}(r)}}^{\rho
,n},\rho _{s,\eta _{n}(r)}(\mu ))W_{\mathcal{\rho }}(dv,dz,\red{dr})  \notag \\
&&+\int_{s}^{t}b(X_{s,\eta _{n}(r)}^{\rho ,n},\mathcal{\rho }_{s,\eta
_{n}(r)}(\mu ))\red{dr}.  \notag
\end{eqnarray}

Notice that we used the time discretization in the coefficients of the
equation but we keep the Poisson measure $\widetilde{N}_{\mathcal{\rho }}$
and the martingale measure $W_{\rho }$ associated to the initial $\rho $ (as in the initial equation~\eqref{stab1}).

\begin{lemma}
\label{L}Suppose that hypotheses (\ref{a1}), (\ref{a1'}), (\ref{a2}) and (%
\ref{a2'}) hold, and $t\rightarrow \rho _{s,t}(\mu )$ is continuous
for every $t\geq s\geq 0$ and every $\mu \in \mathcal{P}_{2}(H).$ Then 
\begin{equation}
\lim_{n\to \infty} \E \vert X_{s,t}^{\rho }-X_{s,t}^{\rho ,n}\vert_H^{2}=0.
\label{stab3}
\end{equation}
\end{lemma}
{\bf Proof.} The proof is rather standard, so we only give a sketch of it. From~(\ref{stab1}) and~(\ref{stab2}), $X_{s,t}^{\rho }-X_{s,t}^{\rho ,n}$ is a sum of three terms: using that $(a+b+c)^2\le 3(a^2+b^2+c^2)$, the isometry property and then (\ref{a1}), (\ref{a1'}) and (\ref{a2}), we get
\begin{align*}
  &\E \vert X_{s,t}^{\rho }-X_{s,t}^{\rho ,n}\vert_H^{2} \le C \int_s^t \left( \E \vert X_{s,r}^{\rho }-X_{s,\eta_n(r)}^{\rho ,n}\vert_H^{2}  + W_2^2(\rho_{s,r}(\mu),\rho_{s,\eta_n(r)}(\mu) )  \right)dr \\
  &\le C \int_s^t \left( 2 \E \vert X_{s,r}^{\rho }-X_{s,r}^{\rho ,n}\vert_H^{2}  +2\E \vert X_{s,r}^{\rho ,n}- X_{s,\eta_n(r)}^{\rho ,n} \vert_H^{2}  + W_2^2(\rho_{s,r}(\mu),\rho_{s,\eta_n(r)}(\mu) )  \right)dr.
\end{align*}
\red{The sublinear growth properties~(\ref{sublin_b}) and (\ref{sublin_csig}) gives} that $\sup_{s\le r\le t}\E[|X_{s,r}^{\rho ,n}|_H^2]<\infty$ for any $t>s$.
 We show then easily that $\sup_{s<r<t}\E \vert X_{s,r}^{\rho ,n}- X_{s,\eta_n(r)}^{\rho ,n} \vert_H^{2} \to_{n\to \infty} 0$, for any $t>s$. Also, we have  $\sup_{s<r<t}W_2^2(\rho_{s,r}(\mu),\rho_{s,\eta_n(r)}(\mu) )$ by the continuity of $t\red{\mapsto} \rho _{s,t}(\mu )$ and Heine theorem. We conclude by Gronwall lemma. $\square$.

\begin{lemma}
\label{STABILITY}Suppose that hypotheses (\ref{a1}), (\ref{a1'}), (\ref{a2}) and (%
\ref{a2'}) hold. Let $\rho _{s,t}$ and $\rho _{s,t}^{\prime }$, $s<t$ be two families of maps $\mathcal{P}_2(H)\to  \mathcal{P}_2(H)$ such that $t\red{\mapsto} \rho_{s,t}(\mu )$ and $t\red{\mapsto} \rho'_{s,t}(\mu )$ are measurable
for  every $\mu \in \mathcal{P}_{2}(H)$. Then, there
exists a constant $C> 0$ such that for every fixed $n\in \N$%
\begin{eqnarray}
&&W_{2}^{2}(\mathcal{L}(X_{s,t}^{\rho ,n}),\mathcal{L}(X_{s,t}^{\rho
^{\prime },n}))  \label{stab4} \\
&\leq &Ce^{C(t-s)}\int_{s}^{t} \left(W_{2}^{2}(\mathcal{\rho }_{s,\eta
_{n}(r)}(\mu ),\mathcal{\rho }_{s,\eta _{n}(r)}^{\prime }(\mu ))+W_{2}^{2}(%
\mathcal{\rho }_{s,r}(\mu ),\mathcal{\rho }_{s,r}^{\prime }(\mu ))\right)dr \notag
\end{eqnarray}%
The constant $C$ depends on the Lipschitz constant of the coefficients, but it does not depend on $n$.
\end{lemma}


\textbf{Proof.} The proof is quite similar to the one of Lemma~\ref{L}, but there is one main additional difficulty: in~Lemma~\ref{L}, the processes are defined with the same noise while here, $X_{s,t}^{\rho ,n}$ and $X_{s,t}^{\rho' ,n}$ are defined with different Poisson random measures and martingale measures. We then first construct two processes  $x_{s,t}^{\rho ,n}$ and $x_{s,t}^{\rho' ,n}$ that have the same law, but are obtained by using the same driving noises $N^{can}$ and $W^{can}$. 

For every $s<r$ we consider $\Pi_{s,r}(dv_{1},dv_{2})$ the
optimal coupling of $\rho _{s,r}(\mu )$ and $\rho _{s,r}^{\prime }(\mu )$. By~\cite[Corollary 5.22]{Villani}, we can assume that $r\mapsto \Pi_{s,r}$ is measurable. Then, there exists  $\tau
_{s,r}:(0,1)\rightarrow $ $H\times H$ such that 
\begin{equation*}
\int_{H\times H}\varphi (v_{1},v_{2})\Pi
_{s,r}(dv_{1},dv_{2})=\int_{0}^{1}\varphi (\tau _{s,r}(w))dw,
\end{equation*}%
and such that $(r,w)\mapsto \tau _{s,r}(w)$ is measurable (see e.g. Blackwell and Dubins~\cite{[BD]}).


Let $x_{s,t}^{\rho ,n}(\mu )$ be the solution of the
equation 
\begin{eqnarray}
x_{s,t}^{\rho ,n} &=&X+\int_{s}^{t}\int_{(0,1)\times E}c(\tau
_{s,r}^{1}(w),z,x_{s,\eta _{n}(r)-}^{\rho ,n},\mathcal{\rho }_{s,\eta
_{n}(r)}(\mu ))\widetilde{N}^{can}(dw,dz,dr)  \notag\\
&&+\int_{s}^{t}\int_{(0,1)\times E}\sigma (\tau _{s,r}^{1}(w),z,x_{s,\eta
_{n}(r)-}^{\rho ,n},\rho _{s,\eta _{n}(r)}(\mu ))W^{can}(dv,dz,ds) \label{stab5}  \\
&&+\int_{s}^{t}b(x_{s,\eta _{n}(r)}^{\rho ,n},\mathcal{\rho }_{s,\eta
_{n}(r)}(\mu ))ds.  \notag
\end{eqnarray}%
By Lemma~\ref{lem_dist}, we show by induction on $k$ that $x_{s,t}^{\rho ,n}$ has the same distribution that $X_{s,t}^{\rho
,n}(\mu )$ for $t\in[s,\red{t^n_k}]$, \red{$t^n_k=s+\frac kn$}. Therefore, $\mathcal{L}(x_{s,t}^{\rho ,n})= \mathcal{L}(X_{s,t}^{\rho
,n}(\mu ))$ for all $t\ge s$.
 Similarly,  $X_{s,t}^{\rho ^{\prime },n}(\mu )$ has the same law as $%
x_{s,t}^{\rho ^{\prime },n}(\mu ),$ solution of the equation 
\begin{eqnarray}
x_{s,t}^{\rho ^{\prime },n} &=&X+\int_{s}^{t}\int_{(0,1)\times E}c(\tau
_{s,r}^{2}(w),z,x_{s,\eta _{n}(r)-}^{\rho ^{\prime },n},\mathcal{\rho }%
_{s,\eta _{n}(r)}^{\prime }(\mu ))\red{\widetilde{N}^{can}}(dw,dz,dr)  \notag \\
&&+\int_{s}^{t}\int_{(0,1)\times E}\sigma (\tau _{s,r}^{2}(w),z,x_{s,\eta
_{n}(r)-}^{\rho ^{\prime },n},\rho _{s,\eta _{n}(r)}^{\prime }(\mu
))\red{W^{can}}(dv,dz,dr)  \label{stab6}\\
&&+\int_{s}^{t}b(x_{s,\eta _{n}(r)}^{\rho ^{\prime },n},\mathcal{\rho }%
_{s,\eta _{n}(r)}^{\prime }(\mu ))dr.  \notag
\end{eqnarray}%
It follows that 
\begin{equation*}
W_{2}^{2}(\mathcal{L}(X_{s,t}^{\rho ,n}),\mathcal{L}(X_{s,t}^{\rho ^{\prime
},n}))=W_{2}^{2}(\mathcal{L}(x_{s,t}^{\rho ,n}),\mathcal{L}(x_{s,t}^{\rho
^{\prime },n}))\leq \E \vert x_{s,t}^{\rho ,n}-x_{s,t}^{\rho ^{\prime
},n}\vert_H^{2}.
\end{equation*}%

Now that we have the stochastic representation with the same driving noises, the proof is similar to the one of Lemma~\ref{L}, but we give more details here.

We denote 
\begin{eqnarray*}
c_{1}(w,z,r) &=&c(\tau _{s,r}^{1}(w),z,x_{s,\red{\eta _{n}(r)}}^{\rho ,n},%
\mathcal{\rho }_{s,\eta _{n}(r)}(\mu )) \\
c_{2}(w,z,r) &=&c(\tau _{s,r}^{2}(w),z,x_{s,\red{\eta _{n}(r)}}^{\rho ^{\prime
},n},\mathcal{\rho }_{s,\eta _{n}(r)}^{\prime }(\mu )) \\
I_{s,t}^{n} &=&\E\left(\int_{s}^{t}\int_{(0,1)\times E}\left\vert
c_{1}(w,z,r)-c_{2}(w,z,r)\right\vert _H^{2}dw\nu (dz)dr\right),
\end{eqnarray*}%
and%
\begin{eqnarray*}
\sigma _{1}(w,z,r) &=&\sigma (\tau _{s,r}^{1}(w),z,x_{s,\red{\eta _{n}(r)}}^{\rho
,n},\mathcal{\rho }_{s,\eta _{n}(r)}(\mu )) \\
\sigma _{2}(w,z,r) &=&\sigma (\tau _{s,r}^{2}(w),z,x_{s,\red{\eta _{n}(r)}}^{\rho
^{\prime },n},\mathcal{\rho }_{s,\eta _{n}(r)}^{\prime }(\mu )) \\
J_{s,t}^{n} &=&\E\left(\int_{s}^{t}\int_{(0,1)\times E}\left\vert \sigma
_{1}(w,z,r)-\sigma _{2}(w,z,r)\right\vert _H^{2}dw\nu (dz)dr\right),
\end{eqnarray*}%
Finally, we set
\begin{equation*}
R_{s,t}^{n}=(t-s)\E\left(\int_{s}^{t}\left\vert b(x_{s,\eta _{n}(r)}^{\rho ,n},%
\mathcal{\rho }_{s,\eta _{n}(r)}(\mu ))-b(x_{s,\eta _{n}(r)}^{\rho ^{\prime
},n},\mathcal{\rho }_{s,\eta _{n}(r)}^{\prime }(\mu ))\right\vert _H^{2}dr \right).
\end{equation*}%
With this notation
\begin{equation*}
\E\vert x_{s,t}^{\rho ,n}-x_{s,t}^{\rho ^{\prime },n}\vert_H
^{2}\leq 3I_{s,t}^{n}+3J_{s,t}^{n}+3R_{s,t}^{n}
\end{equation*}%
Using the Lipschitz property (\ref{a1}), (\ref{a1'}) and (\ref{a2}), we get
\begin{eqnarray*}
I_{s,t}^{n} &\leq &C\int_{s}^{t}\int_{(0,1)\times E}C^{2}(z)(\left\vert \tau
_{s,r}^{1}(w)-\tau _{s,r}^{2}(w)\right\vert _H^{2}+ \E \left\vert x_{s,\red{\eta _{n}(r)}}^{\rho ,n}-x_{s,\red{\eta _{n}(r)}}^{\rho ^{\prime },n}\right\vert _H^{2}
\\
&&+W_{2}^{2}(\mathcal{\rho }_{s,\eta _{n}(r)}(\mu ),\mathcal{\rho }_{s,\eta
_{n}(r)}^{\prime }(\mu )))dw\nu (dz)dr \\
&\leq &C\int_{s}^{t}(\int_{(0,1)}\vert \tau _{s,r}^{1}(w)-\tau
_{s,r}^{2}(w)\vert_H^{2}dw+\E\vert x_{s,\red{\eta _{n}(r)}}^{\rho
,n}-x_{s,\red{\eta _{n}(r)}}^{\rho ^{\prime },n}\vert_H^{2} \\
&&+W_{2}^{2}(\mathcal{\rho }_{s,\eta _{n}(r)}(\mu ),\mathcal{\rho }_{s,\eta
_{n}(r)}^{\prime }(\mu )))dr
\end{eqnarray*}%
Recall that
\begin{equation*}
\int_{(0,1)}\left\vert \tau _{s,r}^{1}(w)-\tau _{s,r}^{2}(w)\right\vert
_H^{2}dw=W_{2}^{2}(\mathcal{\rho }_{s,r}(\mu ),\mathcal{\rho }_{s,r}^{\prime
}(\mu )), 
\end{equation*}%
so that%
\begin{eqnarray*}
I_{s,t}^{n} &\leq &C\int_{s}^{t}\E\vert x_{s,\red{\eta _{n}(r)}}^{\rho
,n}-x_{s,\red{\eta _{n}(r)}}^{\rho ^{\prime },n}\vert_H^{2} dr \\
&&+C\int_{s}^{t}W_{2}^{2}(\mathcal{\rho }_{s,\eta _{n}(r)}(\mu ),\mathcal{%
\rho }_{s,\eta _{n}(r)}^{\prime }(\mu )))+W_{2}^{2}(\mathcal{\rho }%
_{s,r}(\mu ),\mathcal{\rho }_{s,r}^{\prime }(\mu )))dr.
\end{eqnarray*}%
Similar upper bounds hold for $J_{s,t}^{n}$ and $R_{s,t}^{n}$ so finally we
get%
\begin{eqnarray*}
\E \vert x_{s,t}^{\rho ,n}-x_{s,t}^{\rho ^{\prime },n}\vert_H^{2}
&\leq &C\int_{s}^{t} \E \vert x_{s,\red{\eta _{n}(r)}}^{\rho ,n}-x_{s,\red{\eta _{n}(r)}}^{\rho ^{\prime },n}\vert_H^{2}dr \\
&&+C\int_{s}^{t}W_{2}^{2}(\mathcal{\rho }_{s,\eta _{n}(r)}(\mu ),\mathcal{%
\rho }_{s,\eta _{n}(r)}^{\prime }(\mu ))+W_{2}^{2}(\mathcal{\rho }%
_{s,r}(\mu ),\mathcal{\rho }_{s,r}^{\prime }(\mu ))dr.
\end{eqnarray*}%
By Gronwall\footnote{The proof of Gronwall lemma that consists in iterating the inequality works the same with $\eta_{n}(r)$ instead of~$r$ as follows. From $\psi(t)\le C \int_s^t \psi(\eta_n(r)) dr + \varphi(t)$ with $C,\varphi,\psi\ge 0$ and $\varphi$ nondecreasing, we get $\psi(t)\le C \int_s^t \int_0^{\eta_n(r)}\left( \psi(\eta_n(u))du +\varphi(u) \right)du dr + \varphi(t)\le ...\le Ce^{C(t-s)} \varphi(t) $ , using that $\int_s^{\eta_n(r)}\left( \psi(\eta_n(u))du +\varphi(u) \right)du\le \int_s^{r}\left( \psi(\eta_n(u))du +\varphi(u) \right)du$.} lemma, we obtain
\begin{equation*}
  \E \vert x_{s,t}^{\rho ,n}-x_{s,t}^{\rho ^{\prime },n}\vert_H^{2} \leq Ce^{C(t-s)}\int_{s}^{t}(W_{2}^{2}(\mathcal{\rho }_{s,\eta
_{n}(r)}(\mu ),\mathcal{\rho }_{s,\eta _{n}(r)}^{\prime }(\mu )))+W_{2}^{2}(%
\mathcal{\rho }_{s,r}(\mu ),\mathcal{\rho }_{s,r}^{\prime }(\mu ))dr.
\end{equation*}%
$\square $\bigskip

\subsection{Step 3: proof of Theorem \ref{jump}}\label{Step3}

\textbf {\red{Existence.}} We start by showing the existence of a solution of Equation~(\ref{eq1}). We proceed
as follows: let $\theta$ be the flow  given by Theorem~\ref{existence}, and
consider the equation 
\begin{eqnarray}
X_{s,t}^{\theta }(X) &=&X+\int_{s}^{t}\int_{H\times E}c(v,z,X_{s,r-}^{\theta
},\theta _{s,r}(\mu ))\widetilde{N}_{\theta _{s,t}(\mu )}(dv,dz,dr)
\label{V14} \\
&&+\int_{s}^{t}\int_{H\times E}\sigma (v,z,X_{s,r}^{\theta },\theta
_{s,r}(\mu ))W_{\theta _{s,t}(\mu )}(dv,dz,dr)  \notag \\
&&+\int_{s}^{t}b(X_{s,r}^{\theta },\theta _{s,r}(\mu ))ds  \notag
\end{eqnarray}%
with $\mu =\mathcal{L}(X).$ Our aim is to prove that $\theta _{s,t}(\mu
)=\mathcal{L}(X_{s,t}^{\theta })$, which precisely means that $X_{s,t}^{\theta }(X)$ is the solution of Equation (\ref{eq1}).

We fix $n \in \N$ and we put $\red{t^n_{k}}=s+\frac{k}{n}.$ For $\red{t^n_{k}}\leq t<\red{t^n_{k+1}}$
we define $X_{s,t}^{\theta ,n}$ to be the solution of the equation (\ref%
{stab2}) with $\rho =\theta .$ We recall that by Lemma \ref{L}\ we know that 
\begin{equation}
\lim_{n} \E \vert X_{s,t}^{\theta }-X_{s,t}^{\theta ,n}\vert_H^{2}=0.
\label{e1}
\end{equation}%
We will use Lemma~\ref{STABILITY} with $\rho=\theta $ (given before) and $\rho'=\Theta^n_{s,t}$ where, for $t\ge s$ such that $%
\red{t^n_{k}}\leq t<\red{t^n_{k+1}}$ we define 
\begin{equation*}
\Theta _{s,t}^{n}=\Theta _{\red{t^n_{k-1},t^n_{k}}}\circ ...\circ \Theta _{s,\red{t^n_{1}}} \ (\text{convention }\Theta _{s,t}^{n}=\mathrm{Id} \text{ if } s\le t<\red{t^n_1}),
\end{equation*}%
with $\Theta _{s,t}$ given in (\ref{Euler}). Then we associate the equation 
\begin{eqnarray*}
X_{s,t}^{\Theta ,n}(X) &=&X+\int_{s}^{t}\int_{H\times E}c(v,z,X_{s,\red{\eta _{n}(r)-}}^{\Theta ,n},\Theta _{s,\eta _{n}(r)}^{n}(\mu ))\widetilde{N}%
_{\Theta _{s,t}^{n}(\mu )}(dv,dz,dr) \\
&&+\int_{s}^{t}\int_{H\times E}\sigma (v,z,X_{s,\red{\eta _{n}(r)}}^{\Theta
,n},\Theta _{s,\eta _{n}(r)}^{n}(\mu ))W_{\Theta _{s,t}^{n}(\mu )}(dv,dz,dr)
\\
&&+\int_{s}^{t}b(X_{s,\eta _{n}(r)}^{\Theta ,n},\Theta _{s,\eta
_{n}(r)}^{n}(\mu ))ds
\end{eqnarray*}%
We stress that%
\begin{equation*}
\mathcal{L}(X_{s,t}^{\Theta ,n}(X))=\Theta _{s,t}^{n}(\mu).
\end{equation*}%
Now Lemma \ref{STABILITY} gives%
\begin{eqnarray*}
&&W_{2}^{2}(\mathcal{L}(X_{s,t}^{\theta ,n}),\Theta _{s,t}^{n})=W_{2}^{2}(%
\mathcal{L}(X_{s,t}^{\theta ,n}),\mathcal{L}(X_{s,t}^{\Theta ,n})) \\
&\leq &Ce^{C(t-s)}\int_{s}^{t}\left(W_{2}^{2}(\theta _{s,\eta _{n}(r)}(\mu
),\Theta _{s,\eta _{n}(r)}^{n}(\mu ))+W_{2}^{2}(\theta _{s,r}(\mu ),\Theta
_{s,r}^{n}(\mu ))\right)dr.
\end{eqnarray*}%
By using (\ref{m1}), we obtain on the other hand that \begin{equation}\label{cv_unif_Theta_n}
  \lim_{n\to \infty} \sup_{r\in[s,t]}W_{2}^{2}(\theta_{s,r},\Theta _{s,r}^{n})=0.
\end{equation}
  This gives $\lim_{n\to \infty} W_{2}^{2}(\mathcal{L}(X_{s,t}^{\theta ,n}),\Theta _{s,t}^{n})=0$ from the last inequality. Combined this with~(\ref{e1}) and again~(\ref{cv_unif_Theta_n}), we get $\lim_{n\to \infty}  W_{2}^{2}(\mathcal{%
\theta }_{s,t},\Theta _{s,t}^{n})=0$ and thus $\mathcal{L}%
(X_{s,t}^{\theta })=\theta _{s,t}$. This shows the existence of a solution to~(\ref{eq1}).

\smallskip

\textbf{Uniqueness.}
Let $(\overline{X}_{s,t})_{t\ge s}$ be a solution of~(\ref{eq1}) such that $\bar{\theta}_{s,t}(\mu)=\mathcal{L}(\overline{X}_{s,t})$ is such that $t\mapsto \bar{\theta}_{s,t}(\mu)$ is continuous.  From Lemma~\ref{L}, we have 
\begin{eqnarray*}
&&W_{2}^{2}(\mathcal{\theta }_{s,t}(\mu),\overline{\theta }_{s,t}(\mu))=W_{2}^{2}(%
\mathcal{L}(X_{s,t}^{\theta }),\mathcal{L}(X_{s,t}^{\overline{\theta }})) =\lim_{n}W_{2}^{2}(\mathcal{L}(X_{s,t}^{\theta ,n}),\mathcal{L}(X_{s,t}^{%
\overline{\theta },n})).\end{eqnarray*}
Using then Lemma~\ref{STABILITY}, we get
\begin{eqnarray*}
  W_{2}^{2}(\mathcal{L}(X_{s,t}^{\theta ,n}),\mathcal{L}(X_{s,t}^{%
  \overline{\theta },n}))&&\leq Ce^{C(t-s)}\int_{s}^{t}(W_{2}^{2}(\mathcal{\theta }_{s,\eta
_{n}(r)}(\mu ),\overline{\theta }_{s,\eta _{n}(r)}(\mu )))+W_{2}^{2}(%
\mathcal{\theta }_{s,r}(\mu ),\overline{\theta }_{s,r}(\mu ))dr \\
\underset{n\to \infty}{\longrightarrow}&&2Ce^{C(t-s)}\int_{s}^{t}W_{2}^{2}(\mathcal{\theta }_{s,r}(\mu ),\overline{%
\theta }_{s,r}(\mu ))dr.
\end{eqnarray*}
This gives  $W_{2}^{2}(\mathcal{\theta }_{s,t}(\mu),\overline{\theta }_{s,t}(\mu))\le 2Ce^{C(t-s)}\int_{s}^{t}W_{2}^{2}(\mathcal{\theta }_{s,r}(\mu ),\overline{%
\theta }_{s,r}(\mu ))dr$, and therefore, by  Gronwall lemma, $W_{2}^{2}(\mathcal{\theta }_{s,t}(\mu),\overline{\theta }_{s,t}(\mu))=0$.

\subsection{Examples}

In this section we give two simple examples of Equation~\ref{eq1}.

We consider first the case $H=\R^d$ and coefficients
\begin{equation*}
c,\sigma : \R^{d}\times E\times \R^{d}\times \mathcal{P}_{2}(\R^{d})\rightarrow
\R^{d},\quad b: \R^{d}\times \mathcal{P}_{2}(\R^{d})\rightarrow \R^{d}
\end{equation*}%
that satisfy the Lipschitz conditions (\ref{a1})--(\ref{a2'}). Equation~(\ref{eq1}) is in this framework
\begin{eqnarray*}
X_{s,t} &=&X+\int_{s}^{t}\int_{\R^{d}\times E}c(v,z,X_{s,r-},\mathcal{L}%
(X_{s,r}))\widetilde{N}_{\mathcal{L}(X_{s,r})}(dv,dz,dr)  \label{ex1} \\
&&+\int_{s}^{t}\int_{\R^{d}\times E}\sigma (v,z,X_{s,r-},\mathcal{L}%
(X_{s,r}))W_{\mathcal{L}(X_{s,r})}(dv,dz,ds)  \notag \\
&&+\int_{s}^{t}b(X_{s,r},\mathcal{L}(X_{s,r}))ds  \notag
\end{eqnarray*}%
where the law of the starting condition $X$ is $\mu \in \mathcal{P}%
_{2}(\R^{d})$. The same example was considered in~\cite{[AB]}, with $\sigma=0$ and  $N_{\mathcal{L}(X_{s,r})}$ instead of $\widetilde{N}_{\mathcal{L}(X_{s,r})}$, with the function $C$ in (\ref{a1}) satisfying $%
\int_{E}C(z)d\nu (z)<\infty$ instead of $\int_{E}C^2(z)d\nu (z)<\infty$. Then, the process has  finite variation, and~\cite{[AB]} uses a ``standard" sewing lemma to get the existence of a flow. In the present framework, we need instead to use the stochastic sewing lemma to deal with the martingale increments. Thus, Theorem~\ref{jump} gives us the existence of a solution to the above equation, as well as the existence of a flow for the marginal distributions.

We now present an example with an Hilbert space with infinite dimension. We consider the Hilbert space 
\begin{equation*}
H=L^{2}(\rho )=\left\{h:\R^{d}\rightarrow \R^{d}:\int_{\R^{d}}\left\vert
h(x)\right\vert ^{2}\rho (dx)<\infty \right\},
\end{equation*}
where $\rho \in \mathcal{P}_2(\R^d)$ is given. The scalar product on $H$ is $\langle h,h'\rangle_H= \int_{\R^{d}}
  h(x)\red{h'(x)}\rho (dx)$
Then, we consider functions $\overline{c},\overline{\sigma}: \R\times E\times \R^{d}\times
\R\rightarrow \R^{ d}$ and $\bar{b}:\R\times\R \to \R$, and we define the function $c,\sigma:  H\times E\times H\times \mathcal{P}_{2}(H)\rightarrow H$ and $b:H\times  \mathcal{P}_{2}(H)\rightarrow H$
by%
\begin{align}\label{ex2}
c(v,z,h,\mu )(x)&=\overline{c}\left(\left\langle v^c_0 ,v\right\rangle
_{H},z,h(x),\int_{H}\left\langle w^c_0 ,w\right\rangle _{H}\mu (dw) \right),\\
\sigma(v,z,h,\mu )(x)&=\overline{\sigma}\left(\left\langle v^\sigma_0 ,v\right\rangle
_{H},z,h(x),\int_{H}\left\langle w^\sigma_0 ,w\right\rangle _{H}\mu (dw) \right) \notag\\
b(h,\mu)(x) &= \bar{b}\left(h(x), \int_{H}\left\langle w^b_0 ,w\right\rangle _{H}\mu (dw)\right)\notag
\end{align}%
where $v^c_0,v^\sigma_0,w^c_0,w^\sigma_0,w^b_0\in H$ are given.  We consider the equation 
\begin{align}
X_{s,t}=&X+\int_{s}^{t}b(X_{s,r-},\mathcal{L}(X_{s,r})) dr+\int_{s}^{t}\int_{H\times E}\sigma(v,z,X_{s,r-},\mathcal{L}(X_{s,r}))%
W_{\mathcal{L}(X_{s,r})}(dv,dz,dr) \notag \\
&+\int_{s}^{t}\int_{H\times E}c(v,z,X_{s,r-},\mathcal{L}(X_{s,r}))%
\widetilde{N}_{\mathcal{L}(X_{s,r})}(dv,dz,dr)\label{eq_poisson_hilbert}
\end{align}%
where $N_{\mathcal{L}(X_{s,r})}$ (resp. $W_{\mathcal{L}(X_{s,r})}$) is a Poisson point (resp. martingale) measure on $H\times
E\times (s,\infty )$ with intensity%
\begin{equation*}
\mathcal{L}(X_{s,r})(dv)\nu
(dz)dr
\end{equation*}%
and $X\sim \mu \in \mathcal{P}_{2}(H)$. Such type of equations has been studied by Ahmed and Ding~\cite{AhDi} with $c=0$ and a constant function $\sigma$. It has been further studied by Ahmed~\cite{Ahmed} who presents an application to mobile communication. 

We present now the hypotheses that we assume on $\overline{c},\overline{\sigma},\overline{b}$ in
order to ensure that $c,\sigma,b$ verify (\ref{a1})--(\ref{a2'}).  We assume that 
\begin{eqnarray}
\left\vert \overline{c}(a,z,x,\xi )-\overline{c}(a^{\prime },z,x^{\prime
},\xi ^{\prime })\right\vert &\leq &C(z)(\left\vert a-a^{\prime }\right\vert
+\left\vert x-x^{\prime }\right\vert +\left\vert \xi -\xi ^{\prime
}\right\vert ), \label{assump_cbar}\\
\text{with}\quad \int_{E}\red{\big(C^{2}(z)+|\bar{c}(0,z,0,\xi_0)|^2\big)}\nu (dz) &<&\infty , \notag
\end{eqnarray}%
where $\xi_0\in \R$ is such that $\xi_0=\int_{H} \langle w^c_0, w \rangle
\mu_0(dw)$ for some $\mu_0 \in \mathcal{P}_2(H)$. 
This gives in particular that $c:H\times E\times H\times \mathcal{P}_{2}(H)\rightarrow H$. We assume that $\overline{\sigma}$ satisfies the same condition as $\overline{c}$ and that $\bar{b}$ is uniformly Lipschitz on $\R \times\R$. 

We first check that $c$ (and thus $\sigma$) satisfies (\ref{a1})--(\ref{a2'}). We have 
\begin{eqnarray}
 &&  \quad \quad  \quad \quad \left\vert c(v,z,h,\mu )(x)-c(v^{\prime },z,h^{\prime },\mu^{\prime
})(x)\right\vert  \label{ex4} \\
&\leq &\left\vert \overline{c}\left(\left\langle v^c_0 ,v\right\rangle
_{H},z,h(x),\int_{H}\left\langle w^c_0 ,w\right\rangle _{H}\mu (dw)\right)-%
\overline{c}\left(\left\langle v^c_0 ,v^{\prime }\right\rangle _{H},z,h^{\prime
}(x),\int_{H}\left\langle w^c_0 ,w\right\rangle _{H}\mu ^{\prime
}(dw)\right)\right\vert  \notag \\
&\leq &C(z) \bigg(\left\vert \left\langle v^c_0 ,v\right\rangle _{H}-\left\langle
v^c_0 ,v^{\prime }\right\rangle _{H}\right\vert +\left\vert h(x)-h^{\prime
}(x)\right\vert  \notag \\
&&+\left\vert \int_{H}\left\langle w^c_0 ,w\right\rangle _{H}\mu
(dw)-\int_{H}\left\langle w^c_0 ,w\right\rangle _{H}\mu ^{\prime
}(dw)\right\vert \bigg).  \notag
\end{eqnarray}%
Taking the square and integrating with respect to~$\rho$, we get
\begin{eqnarray*}
&&\left\vert c(v,z,h,\mu )-c(v^{\prime },z,h^{\prime },\mu ^{\prime
})\right\vert^2_{H} \\
&\leq &3C(z)^2 \bigg(\left\vert \left\langle v_0 ,v\right\rangle _{H}-\left\langle
v_0 ,v^{\prime }\right\rangle _{H}\right\vert^2 +\left\vert h-h^{\prime}\right\vert^2_{H} \\
&&+\left\vert \int_{H} \left\langle w^c_0 ,w \right\rangle _{H}\mu(dw) -\int_{H}\left\langle w^c_0 ,w\right\rangle _{H}\mu ^{\prime}(dw)\right\vert ^2 \bigg).
\end{eqnarray*}%
The first term can be upper bounded by using $\left\vert \left\langle v^c_0 ,v\right\rangle _{H}-\left\langle v^c_0
,v^{\prime }\right\rangle _{H}\right\vert \leq \left\vert v^c_0 \right\vert
_{H}\left\vert v-v^{\prime }\right\vert _{H}$. For the third term, let \red{us} denote $Q(w)=\left\langle
w^c_0 ,w\right\rangle _{H}$ and take $\Pi$ a $W_2$-optimal coupling of $\mu$ and $\mu'$. We have 
\begin{eqnarray*}
  \left\vert \int_{H}\left\langle w^c_0 ,w\right\rangle _{H}\mu
  (dw)-\int_{H}\left\langle w^c_0 ,w\right\rangle _{H}\mu ^{\prime
  }(dw)\right\vert 
  &=&\left\vert \int_{H\times H} \left\langle w^c_0 ,w\right\rangle _{H} - \left\langle w^c_0 ,w'\right\rangle _{H} \Pi(dw,dw') \right\vert\\
  &\le& \left\vert \int_{H\times H} |w^c_0|_H|w-w'|_H  \Pi(dw,dw') \right\vert \\
  &\le&|w^c_0| \left\vert \int_{H\times H} |w-w'|^2_H  \Pi(dw,dw') \right\vert^{1/2} \\
  &=&|w^c_0| W_2(\mu,\mu').
\end{eqnarray*}
We conclude that%
\begin{eqnarray*}
&&\left\vert c(v,z,h,\mu )-c(v^{\prime },z,h^{\prime },\mu ^{\prime
})\right\vert _{H} \\
&\leq &\sqrt{3} C(z)(\left\vert v^c_0 \right\vert _{H}\left\vert v-v^{\prime
}\right\vert _{H}+\left\vert h-h^{\prime }\right\vert _{H}+\left\vert
w^c_0 \right\vert _{H}W_{2} (\mu ,\mu ^{\prime })).
\end{eqnarray*}%
So the Lipschitz property (\ref{a1})--(\ref{a2'}) is verified by $c$ and $\sigma$. We check in the same way the properties for~$b$. By Theorem~\ref{jump}, we get the existence of a solution to~(\ref{eq_poisson_hilbert}), and the uniqueness of the flow $\theta_{s,t}(\mu)$ that describes the marginal laws of this process.

\bigskip

\appendix

\section{ Technical results}

Hereafter, $\kp$ will denote a universal constant depending only on $p$.

\subsection{Proof of Lemma \ref{BASIC1}.}\label{AppendixlemBASIC1}

\smallskip

Let $(X_k,0\le k\le N)$ be defined by~(\ref{def_Xk}). We set $\eta _{i}=\mathcal{E}^1_{t_{i},r_i,t_{i+1}}(X_{i}),\widehat{\eta }_{i}=\E(\eta _{i}|{\mathcal{F}%
_{t_{i}}})$ and $\widetilde{\eta }_{i}=\eta _{i}-\widehat{\eta }%
_{i}.$ Then by {\bf G}, we get%
\begin{equation}\label{G0bis}
\left\Vert \widehat{\eta }_{i}\right\Vert _{p}\leq
L_{sew}(t_{i+1}-t_{i})(1+\left\Vert X^1_{i}\right\Vert _{p})\quad \text{ and }\quad
\left\Vert \widetilde{\eta }_{i}\right\Vert _{p}\leq
L_{sew}(t_{i+1}-t_{i})^{1/2}(1+\left\Vert X_{i}^1\right\Vert _{p}).
\end{equation}
For $k\ge 0$, we rewrite $\red{X^1_{k+1}}$ as follows
\begin{equation*}
X^1_{k+1}=X^1_{0}+\sum_{i=0}^{k}\eta _{i}=\red{X^1_{0}}+\widehat{T}_{k+1}+\widetilde{T}%
_{k+1},
\end{equation*}%
with $\widehat{T}_{k+1}=\sum_{i=0}^{k}\widehat{\eta }_{i}$ and $\widetilde{T}%
_{k+1}=\sum_{i=0}^{k}\widetilde{\eta }_{i}$. Therefore we get 
\begin{equation*}
\left( 1+\left\Vert X^1_{k+1}\right\Vert _{p}\right) ^{p}\leq 3^{p-1}\left(
1+\left\Vert X^1_{0}\right\Vert _{p}\right) ^{p}+3^{p-1}\left\Vert \widehat{T}%
_{k+1}\right\Vert _{p}^{p}+3^{p-1}\left\Vert \widetilde{T}_{k+1}\right\Vert
_{p}^{p}.
\end{equation*}%
By the triangle inequality,~(\ref{G0bis}) and then Jensen inequality, we have 
\begin{align*}
\left\Vert \widehat{T}_{k+1}\right\Vert _{p}^{p}& \leq \left(
\sum_{i=0}^{k}\left\Vert \widehat{\eta }_{i}\right\Vert _{p}\right)
^{p}\leq \left( \sum_{i=0}^{k}L_{sew}(t_{i+1}-t_{i})(1+\left\Vert
X^1_{i}\right\Vert _{p})\right) ^{p} \\
& \leq
L_{sew}^{p}(t_{k+1}-t_{0})^{p-1}\sum_{i=0}^{k}(t_{i+1}-t_{i})(1+\left\Vert
X^1_{i}\right\Vert _{p})^{p}.
\end{align*}%
By using the Burkholder-Davis-Gundy inequality for Hilbert space valued
martingales (see Marinelli and Röckner~\cite[Theorem 1.1]{MR}), 
Jensen inequality and then~(\ref{G0bis}), we get 
\begin{align*}
\Vert \widetilde{T}_{k+1}\Vert _{p}^{p}= \E \Big(\Big|\sum_{i=0}^{k}\widetilde{%
\eta }_{i}\Big|_{H}^{p}\Big)& \leq \kappa _{p}^{p} \E \Big(\Big(\sum_{i=0}^{k}|%
\widetilde{\eta }_{i}|_{H}^{2}\Big)^{p/2}\Big)=\kappa _{p}^{p} \E\Big(\Big(%
\sum_{i=0}^{k}(t_{i+1}-t_{i})\frac{|\widetilde{\eta }_{i}|_{H}^{2}}{%
(t_{i+1}-t_{i})}\Big)^{p/2}\Big) \\
& \leq \kappa _{p}^{p}(t_{k+1}-t_{0})^{\frac{p}{2}-1}%
\sum_{i=0}^{k}(t_{i+1}-t_{i})\frac{\E(|\widetilde{\eta }_{i}|_{H}^{p})}{%
(t_{i+1}-t_{i})^{p/2}} \\
& \leq L_{sew}^{p}\kappa _{p}^{p}(t_{N}-t_{0})^{\frac{p}{2}%
-1}\sum_{i=0}^{k}(t_{i+1}-t_{i})(1+\Vert X^1_{i}\Vert _{p})^{p}
\end{align*}%
%
Therefore, \red{since $t_N-t_0\leq T$,} we get 
\begin{align*}
(1+\left\Vert X^1_{k+1}\right\Vert _{p})^{p}\leq & 3^{p-1}(1+\left\Vert
X^1_{0}\right\Vert _{p})^{p} \\
& +3^{p-1}L_{sew}^{p} \red{T^{\frac{p}{2}-1}}
\left( \red{T^{%
	\frac{p}{2}}}
+\kappa _{p}^{p}\right) \sum_{i=0}^{k}(t_{i+1}-t_{i})(1+\Vert
X^1_{i}\Vert _{p})^{p},
\end{align*}%
and thus 
\begin{align*}
  (1+\left\Vert X^1_{k+1}\right\Vert _{p})^{p}&\leq  3^{p-1}(1+\left\Vert
X^1_{0}\right\Vert _{p})^{p} \exp\left(3^{p-1}L_{sew}^{p}
\red{T^{\frac{p}{2}-1}}
\left( \red{T^{\frac{p}{2}}}
+\kappa _{p}^{p}\right) (t_{k+1}-t_0) \right)  \\
&\leq  3^{p-1} \exp\left(3^{p-1}L_{sew}^{p}T^{\frac{p}{2}}\left( T^{\frac{p}{2}}+\kappa _{p}^{p}\right)  \right)(1+\left\Vert
  X^1_{0}\right\Vert _{p})^{p}
\end{align*}
 by Gronwall lemma. This gives~(\ref{majo_Lambda}) since $X^1_0\sim \mu$ and $X^1_N\sim \Theta^{\pi}_{t_0,t_N}(\mu)$.

\subsection{Proof of Lemma \ref{sewing}.}\label{AppendixlemBASIC2}
We note $\pi=\{t_{0}<t_{1}<....<t_{N}\}$ and $\pi'=\{t_{0}\leq r_{0}<t_{1}<....<t_{N-1}\leq r_{N-1}<t_{N}\}$ the simple subpartition of $\pi$.

Let $\Pi (\mu ,\nu )\in \mathcal{P}_p(H\times H)$
be an optimal $W_{p}$ coupling of $\mu$ and $\nu$, and let $%
X_{0}=(X_{0}^{1},X_{0}^{2})\in L_{t_{0}}^{p}\times L_{t_{0}}^{p}$ with law $%
\Pi (\mu ,\nu )$ (this is possible by hypothesis \textbf{A}). This
means that
\begin{equation*}
\left\Vert X_{0}^{1}-X_{0}^{2}\right\Vert _{p}^{p}=\int_{H\times
H}\left\vert x-y\right\vert _{H}^{p}\Pi (\mu ,\nu )(dx,dy)=W_{p}^{p}(\mu
,\nu ).
\end{equation*}
Let $(X_k,0\le k\le N)$ be defined by~(\ref{def_Xk}) and denote%
\begin{equation*}
\delta_{i} =\mathcal{E}_{t_{i},r_i,t_{i+1}}^{2}(X^2_{i})-\mathcal{E}_{
t_{i},r_i,t_{i+1}}^{1}(X^1_{i})\quad \mbox{and}\quad \widehat{\delta }_{i} =\E_{%
\mathcal{F}_{t_{i}}}(\delta _{i}),\quad \widetilde{\delta }_{i}=\delta _{i}-%
\widehat{\delta }_{i}.
\end{equation*}
Then by (\ref{V7}) and (\ref{V8}), we get%
\begin{eqnarray*}
\left\Vert \widehat{\delta }_{i}\right\Vert _{p} &\leq
&C_{sew}(t_{i+1}-t_{i})
\red{\Big(}
\left\Vert X_{i}^{1}-X_{i}^{2}
\right\Vert
_{p}+(1+\|X_i\|_p)\left\vert \pi \right\vert ^{\beta }
\red{\Big)}, \\
\left\Vert \widetilde{\delta }_{i}\right\Vert _{p} &\leq
&C_{sew}(t_{i+1}-t_{i})^{1/2}
\red{\Big(}\left\Vert X_{i}^{1}-X_{i}^{2}\right\Vert
_{p}+(1+\|X_i\|_p)\left\vert \pi \right\vert ^{\beta}
\red{\Big)}.
\end{eqnarray*}%
We will use the decomposition 
\begin{equation*}
X_{k+1}^{2}-X_{k+1}^{1}=X_{0}^{2}-X_{0}^{1}+\sum_{i=0}^{k}\delta
_{i}=X_{0}^{2}-X_{0}^{1}+\widehat{S}_{k+1}+\widetilde{S}_{k+1},
\end{equation*}%
with%
\begin{equation*}
\widehat{S}_{k+1}=\sum_{i=0}^{k}\widehat{\delta }_{i},\quad \mbox{and}\quad 
\widetilde{S}_{k+1}=\sum_{i=0}^{k}\widetilde{\delta }_{i}.
\end{equation*}

We notice first that by the triangle inequality%
\begin{eqnarray*}
\left\Vert \widehat S_{k+1}\right\Vert _{p}^{p} &\leq &
(\sum_{i=0}^{k}\left\Vert \widehat{\delta }_{i}\right\Vert _{p})^{p} 
\leq
\left(\sum_{i=0}^{k} C_{sew}(t_{i+1}-t_{i})\red{\Big(}\left\Vert
X_{i}^{1}-X_{i}^{2}\right\Vert _{p}+(1+\|X_i\|_p)\left\vert \pi \right\vert
^{\beta } \red{\Big)}\right)^{p} \\
&\leq &2^{p-1}C_{sew}^{p}(t_{k+1}-t_{0})^{p-1}\sum_{i=0}^{k}(t_{i+1}-t_{i})%
\big(\left\Vert X_{i}^{1}-X_{i}^{2}\right\Vert
_{p}^{p}+(1+\|X_i\|_p)^p\left\vert \pi \right\vert ^{p\beta }\big) \\
&\leq
&2^{p-1}C_{sew}^{p}(t_{N}-t_{0})^{p-1}\sum_{i=0}^{k}(t_{i+1}-t_{i})\left\Vert X_{i}^{1}-X_{i}^{2}\right\Vert _{p}^{p} \\
&&+2^{p-1}C_{sew}^{p}(t_{N}-t_{0})^{p}\sup_{i\leq
N}(1+\|X_i\|_p)^p\left\vert \pi \right\vert ^{p\beta }.
\end{eqnarray*}%
Using Lemma~\ref{BASIC1}, \red{that we apply to the restriction of $\pi$ and $\pi'$ on $[t_0,t_i]$,} we have $1+\|X_i\|_p=1+\|X^1_i\|_p+\|X^2_i\|_p\le A_T(2+\|X^1_0\|_p+\|X^2_0\|_p )=A_T(2+\|X_0\|_p )$ and thus
 \begin{eqnarray*}
\left\Vert \widehat{S}_{k+1}\right\Vert _{p}^{p} &\leq
&2^{p-1}C_{sew}^{p}
\red{T^{p-1}}
\sum_{i=0}^{k}(t_{i+1}-t_{i})\left\Vert X_{i}^{1}-X_{i}^{2}\right\Vert _{p}^{p} \\
&&+2^{p-1}C_{sew}^{p}A_T^p (2+\|X_{0}%
\|_p)^{p}
\red{T^p}
\left\vert \pi \right\vert ^{p\beta }.
\end{eqnarray*}

Moreover, using the Burkholder-Davis-Gundy inequality and then Jensen
inequality, we get 
\begin{align*}
\|\widetilde{S}_{k+1}\|_p^{p} &=\E\Big(\Big|\sum_{i=0}^{k}\widetilde{\delta}_i%
\Big|_H^p\Big) \leq \kappa_p^p \E\Big(\Big(\sum_{i=0}^{k}|\widetilde{\delta}%
_i|_H^2\Big)^{p/2}\Big) =\kappa_p^p \E\Big(\Big(\sum_{i=0}^{k}(t_{i+1}-t_i) 
\frac{|\widetilde{\delta}_i|_H^2}{t_{i+1}-t_i}\Big)^{p/2}\Big) \\
&\leq \kappa_p^p(t_{k+1}-t_0)^{\frac p 2 -1}\sum_{i=0}^k(t_{i+1}-t_i)\frac{\E%
\big(|\widetilde{\delta}_i|_H^p\big)}{(t_{i+1}-t_i)^{p/2}} \\
&\leq C_{sew}^p\kappa_p^p (t_{N}-t_0)^{\frac{p-2}{2}}
\sum_{i=0}^k(t_{i+1}-t_i)\big(\left\Vert X_{i}^{1}-X_{i}^{2}\right\Vert
_{p}+(1+\|X_i\|_p)\left\vert \pi \right\vert ^{\beta }\big)^{p}
\end{align*}
Therefore, we get by using again Lemma \ref{BASIC1}
\begin{align*}
\|\widetilde{S}_{k+1}\|_p^{p} &\leq 2^{p-1}C_{sew}^p\kappa_p^p
\red{T^{\frac{p}{2} -1} }
\sum_{i=0}^k(t_{i+1}-t_i)\left\Vert
X_{i}^{1}-X_{i}^{2}\right\Vert _{p}^{p} \\
& \quad + 2^{p-1}C_{sew}^p\kappa_p^p \red{T^{\frac{p}{2} }}
A_T^p (2+\|X_{0}\|_p)^{p} \left\vert \pi \right\vert
^{p\beta }.
\end{align*}
Finally, inserting the above estimates for $\widehat S_{k+1}$ and $%
\widetilde S_{k+1}$, we obtain 
\begin{align*}
\|X^1_{k+1}-X^2_{k+1}\|_p^p &\leq
3^{p-1}\|X^1_{0}-X^2_{0}\|_p^p+3^{p-1}\|\widehat
S_{k+1}\|_p^p+3^{p-1}\|\widetilde S_{k+1}\|_p^p \\
&\leq
3^{p-1}\|X^1_{0}-X^2_{0}\|_p^p+\tilde{a}_T(1+\|X_0\|_p)^p
|\pi|^{p\beta} +\tilde{b}_T \sum_{i=0}^k (t_{i+1}-t_i)\left\Vert
X_{i}^{1}-X_{i}^{2}\right\Vert _{p}^{p}
\end{align*}
where 
\begin{align*}
&\tilde{a}_T =6^{p-1}2^pC^p T^{\frac{p}{2}}(T^{\frac{p}{2}%
}+\kappa_p^p) A_T^p , \ \tilde{b}_T = 6^{p-1}C_{sew}^p T^{\frac{p}{2}-1}\big(T^{\frac p2}+\kappa_p^p\big).
\end{align*}
The Gronwall's Lemma then gives 
\begin{equation*}
\sup_{k\leq N}\|X^1_{k}-X^2_{k}\|_p^p \leq \Big(3^{p-1}\|X^1_{0}-X^2_{0}
\|_p^p+\tilde{a}_T (1+\|X_0\|_p)^p |\pi|^{p\beta}%
\Big)\exp\big(\tilde{b}_{T} T\big)
\end{equation*}
Since $X^1_N\sim \Theta^{\pi}(\mu)$ and $X^2_N\sim\Theta^{\pi'}(\nu)$, we have $W_p(\Theta^{\pi}(\mu),\Theta^{\pi'}(\nu))\le \|X^1_{N}-X^2_{N}\|_p$. Besides, we have $\|X^1_{0}-X^2_{0}\|_p^p=W_p^p(\mu,\nu)$, and we get~\eqref{V9'}. 
$\square $

\subsection{Continuity}

We recall the notation $\eta _{n}(t)=\frac{k}{2^{n}}$ for $\frac{k}{2^{n}}%
\leq t<\frac{k+1}{2^{n}}.$

\begin{lemma}
\label{lemma-app1}  Let $\Theta_{s,t}:\mathcal{P}_p(H)\to \mathcal{P}_p(H)$ be a family of operators satisfying the $\beta$-sewing property and \red{the continuity property} $\mathbf{C}$. Let $t_m<T$ and $\pi =\{t_{0}<...<t_{m}\}$ a partition. For $n$
such that $2^{-n}<\min_{1\le \ell \le m} (t_{\ell}-t_{\ell-1})$, we define
the ``projected partition" $\pi_{n}=\{\eta_{n}(t_{0})<...<\eta
_{n}(t_{m})\} $. Then, there is a
constant~$C_T\in \R_+$ depending on 
$T$, $p$, $L_{sew}$ and $C_{sew}$ such that
\begin{equation*}
W_{p}(\Theta _{t_{0},t_{m}}^{\pi }(\mu ),\Theta _{\eta _{n}(t_{0}),\eta
_{n}(t_{m})}^{\pi _{n}}(\mu ))\leq C_T (1+\|\mu\|_p) \Big(\varphi\Big(\frac
1{2^n}\Big)+\Big(t_m-t_0+\frac 1{2^n}\Big)^{1/2}\Big(|\pi|+\frac 1{2^n}\Big)%
^{\beta}\Big).
\end{equation*}
\end{lemma}

\textbf{Proof } \red{Hereafter, $C_T\in\R_+$ denotes a constant, possibly depending on 
	$T$, $p$, $L_{sew}$ and $C_{sew}$, that may vary from a line to another.} 

We construct the partitions 
\begin{eqnarray*}
\pi _{n}^{\prime } &=&\{\eta _{n}(t_{0})\leq t_{0}<\eta _{n}(t_{1})\leq
t_{1}<....<\eta _{n}(t_{m-1})\leq t_{m-1}<\eta _{n}(t_{m})\} \\
\pi ^{\prime } &=&\{t_{0}<\eta _{n}(t_{1})\leq t_{1}<....<\eta
_{n}(t_{m-1})\leq t_{m-1}<\eta _{n}(t_{m})\leq t_{m}\}
\end{eqnarray*}%
and we notice that $\pi_{n}^{\prime }$ (resp. $\pi^{\prime }$) is a simple
subpartition of $\pi _{n}$ (resp. $\pi$). Therefore, we can apply Lemma \ref%
{sewing} and obtain%
\begin{align}
W_{p}(\Theta _{\eta _{n}(t_{0}),\eta _{n}(t_{m})}^{\pi _{n}}(\mu ),\Theta
_{\eta _{n}(t_{0}),\eta _{n}(t_{m})}^{\pi _{n}^{\prime }}(\mu )) &\leq
C_T (1+2\|\mu\|_p)\Big(t_{m}-t_{0}+\frac{1}{2^{n}}%
\Big)^{1/2}\Big(\left\vert \pi\right\vert +\frac 1{2^n}\Big)^{\beta},
\label{app1}\\
W_{p}(\Theta _{t_{0},t_{m}}^{\pi}(\mu ),\Theta
_{t_{0},t_{m}}^{\pi^{\prime }}(\mu )) &\leq
C_T (1+2\|\mu\|_p) \Big(t_{m}-t_{0}%
\Big)^{1/2}\left\vert \pi\right\vert^{\beta},
\label{app2}
\end{align}%
where we have used that $|\pi_n|\leq |\pi|+2^{-n}$.

We introduce now the partition 
\begin{equation*}
\widehat{\pi }=\{t_{0}<\eta _{n}(t_{1})\red{\leq}t_{1}<....<\eta
_{n}(t_{m-1})\red{\leq}t_{m-1}<\eta _{n}(t_{m})\}
\end{equation*}%
and we notice 
\begin{eqnarray*}
\Theta _{\eta _{n}(t_{0}),\eta _{n}(t_{m})}^{\pi _{n}^{\prime }}(\mu )
&=&\Theta _{t_{0},\eta _{n}(t_{m})}^{\widehat{\pi }} \circ \Theta _{\eta
_{n}(t_{0}),t_{0}}(\mu ) \\
\Theta _{t_{0},t_{m}}^{\pi ^{\prime }}(\mu ) &=&\Theta _{\eta
_{n}(t_{m}),t_{m}} \circ \Theta _{t_{0},\eta _{n}(t_{m})}^{\widehat{\pi }}(\mu ).
\end{eqnarray*}
By using the asymptotic Lipschitz property (\ref{V9"}), the continuity 
$\mathbf{C}$ and Lemma~\ref{BASIC1}, we also easily obtain 
\begin{eqnarray*}
&&W_{p}(\Theta _{\eta _{n}(t_{0}),\eta _{n}(t_{m})}^{\pi _{n}^{\prime }}(\mu
),\Theta _{t_{0},\eta _{n}(t_{m})}^{\widehat{\pi }}(\mu )) \\
&&\leq C_T W_p(\Theta _{\eta
_{n}(t_{0}),t_{0}}(\mu ),\mu )+C_T(1+\|\mu\|_p+\|\Theta _{\eta
_{n}(t_{0}),t_{0}}(\mu ) \|_p)(\eta_n(t_m)-t_0)^{1/2}|\pi|^\beta \\
&&\leq C_T \Big((1+\|\mu\|_p) \varphi\Big(\frac 1{2^n}\Big)+ (1+\|\mu\|_p+\|\Theta _{\eta
_{n}(t_{0}),t_{0}}(\mu ) \|_p) (t_m-t_0)^{1/2}|%
\pi|^{\beta}\Big) \\
&&\le C_T (1+\|\mu\|_p)\Big( \varphi\Big(\frac 1{2^n}\Big) +(A_T+1) (t_m-t_0)^{1/2}|%
\pi|^{\beta}\Big) \\
&&W_{p}(\Theta _{t_{0},t_{m}}^{\pi ^{\prime }}(\mu ),\Theta _{t_{0},\eta
_{n}(t_{m})}^{\widehat{\pi }}(\mu )) \leq  A_T (1+\|\mu\|_p)\varphi\Big(\frac 1{2^n}\Big)
\end{eqnarray*}%
 Then, gathering these estimates with~(\ref{app1}) and~(\ref%
{app2}) in the triangular inequality 
\begin{align*}
W_{p}(&\Theta_{t_{0},t_{m}}^{\pi }(\mu ),\Theta _{\eta _{n}(t_{0}),\eta
_{n}(t_{m})}^{\pi _{n}}(\mu )) \le W_{p}(\Theta _{t_{0},t_{m}}^{\pi }(\mu
),\Theta _{t_{0},t_{m}}^{\pi^{\prime }}(\mu ))+ W_{p}(\Theta
_{t_{0},t_{m}}^{\pi^{\prime }}(\mu ),\Theta _{t_{0},\eta _{n}(t_{m})}^{%
\widehat{\pi}}(\mu )) \\
&+W_{p}(\Theta _{t_{0},\eta _{n}(t_{m})}^{\widehat{\pi}}(\mu ), \Theta
_{\eta _{n}(t_{0}),\eta _{n}(t_{m})}^{\pi^{\prime }_{n}}(\mu ))+W_{p}(
\Theta _{\eta _{n}(t_{0}),\eta _{n}(t_{m})}^{\pi^{\prime }_{n}}(\mu ) ,
\Theta _{\eta _{n}(t_{0}),\eta _{n}(t_{m})}^{\pi _{n}}(\mu )),
\end{align*}
we get the claim. $\square $

\subsection{Equalities in distribution for Euler schemes}\label{append_schemes}
 
\begin{lemma}\label{lem_dist}Let $c,\sigma:H\times E \times H \times \mathcal{P}_2(H)$ be coefficients satisfying the sublinear property~\eqref{sublin_csig}. Let $\eta _{s_0,t}\in \mathcal{P}_{2}(H)$, $0\le s_0<t$,  such that $t\mapsto \eta_{s_0,t}$ is measurable, and let $N_{\mathcal{\eta}}$ (resp. $W_{\mathcal{\eta}}$) be a Poisson point (resp. martingale) measure on $%
 H\times E\times (s_0,\infty )$ with intensity
 $\eta_{s_0,r}(dv)\nu (dz)1_{(s_0,\infty
 )}(r)dr$.
 Let  $N^{can}$ (resp. $W^{can}$)   be a Poisson point (resp. martingale) measure on   $(0,1)\times E \times (0,\infty)$ with intensity $1_{(0,1)(w)}dw \nu(dz)1_{(0,\infty
 )}(r)dr$.   Let $\tau_{s_0,u}:(0,1)\to H$ be such that $\tau_{s_0,\red{u}}(U)\sim \eta_{s_0,\red{u}}$ for $U$ uniformly distributed on $(0,1)$.
  
  Then, for any $x\in H$, $\rho \in \mathcal{P}_2(H)$ and $s_0\le s<t$,  
  \begin{align*}
    &\int_s^t \int_{H\times E} c(v,z,x,\rho) \red{\tilde N_{\eta}}(dv,dz,du)\sim \int_s^t \int_{H\times E} c(\tau_{s_0,u}(w),z,x,\rho) \red{\tilde N^{can}}(dw,dz,du), \\
    &\int_s^t \int_{H\times E} \sigma(v,z,x,\rho) W_{\eta}(dv,dz,du)\sim \int_s^t \int_{H\times E} \sigma(\tau_{s_0,u}(w),z,x,\rho) W^{can}(dw,dz,du).
  \end{align*}
\end{lemma}
{\bf Proof.} Let us consider $h \in H$, $
  I_{s,t}=\int_s^t \int_{H\times E} c(v,z,x,\rho) \red{\tilde N_{\eta}}(dv,dz,du)$ and $I^{can}_{s,t}=\int_s^t \int_{H\times E} c(\tau_{s_0,u}(w),z,x,\rho) \red{\tilde N^{can}}(dw,dz,du).$ 
 We show the distribution equality by calculating the characteristic functions. We apply Itô's formula~\cite[Theorem 5.1 p. 67]{[IW]} (see~\cite{[AZ]} on Hilbert space) to $e^{i\langle h,I_{s,t}\rangle}$. \red{To simplify the notation, we simply write $\<\cdot,\cdot\>$ in place of $\<\cdot,\cdot\>_H$.  By taking the expectation we get }
\begin{align*}
  \E[e^{i\langle h,I_{s,t}\rangle}]&=1 +\E\left[ \int_s^t \int_{H\times E} \left(e^{i\langle h,I_{s,u}+c(v,z,x,\rho)\rangle} - e^{i\langle h,I_{s,u}\rangle} -i\red{\langle  h, c(v,z,x,\rho)\rangle} e^{i\langle h,I_{s,u}\rangle}\right) \eta_{s_0,u}(dv)\nu(dz)du \right] \\
  &=1+\int_s^t \E[ e^{i\langle h,I_{s,u}\rangle}] \left(  \int_{H\times E} \left(e^{i\langle h,c(v,z,x,\rho)\rangle} - 1 -i\red{\langle h, c(v,z,x,\rho)\rangle} \right) \eta_{s_0,u}(dv)\nu(dz) \right)du
\end{align*}
Similarly, 
\begin{align*}
  \E[e^{i\langle h,I^{can}_{s,t}\rangle}]&=1+\int_s^t \E[ e^{i\langle h,I^{can}_{s,u}\rangle}] \left(  \int_{(0,1)\times E} \left(e^{i\langle h,c(\tau_{s_0,u}(w),z,x,\rho)\rangle} - 1 -i\red{\langle  h, c(\tau_{s_0,u}(w),z,x,\rho)\rangle}\right) dw \nu(dz) \right)du\\
  &=1+\int_s^t \E[ e^{i\langle h,I^{can}_{s,u}\rangle}] \left(  \int_{H\times E} \left(e^{i\langle h,c(v,z,x,\rho)\rangle} - 1 -i\red{\langle  h, c(v,z,x,\rho)\rangle} \right) \eta_{s_0,u}(dv)\nu(dz) \right)du,
\end{align*}
by a change of variable, using that the distribution of $\tau_{s_0,u}(\red{U})$ with $\red{U}\sim \mathcal{U}([0,1])$ is $\eta_{s_0,u}$. Thus $t\mapsto \E[e^{i\langle h,I_{s,t}\rangle}]$ and $t\mapsto \E[e^{i\langle h,I^{can}_{s,t}\rangle}]$ solve the same ODE for $t\ge s$ and are equal. This shows the first identity. 

For the second one, we define   $J_{s,t}=\int_s^t \int_{H\times E} \sigma(v,z,x,\rho) W_{\eta}(dv,dz,du)$ for $t\ge s$ and $J_{s,t}^{can}=\int_s^t \int_{H\times E} \sigma(\tau_{s_0,u}(w),z,x,\rho) W^{can}(dw,dz,du)$. We 
use~\cite[Corollary I-7]{ElKM} to get
\begin{align*}
  \E[e^{i\langle h,J_{s,t}\rangle \red{+\frac 12} \int_s^t \int_{H\times E}\langle h,\sigma(v,z,x,\rho) \rangle^2 \eta_{s_0,u}(dv)\nu(dz) du } ]=1
\end{align*}
Therefore, 
\begin{align*}
  \E[e^{i\langle h,J_{s,t}\rangle}]= \red{e^{-\frac 12 \int_s^t \int_{H\times E}\langle h,\sigma(v,z,x,\rho) \rangle^2 \eta_{s_0,u}(dv)\nu(dz) du }} ,
\end{align*}
and we get in the same way
\begin{align*} 
  \E[e^{i\langle h,\red{J^{can}_{s,t}}\rangle}]=\red{e^{-\frac 12 \int_s^t \int_{(0,1)\times E}\langle h,\sigma(\tau_{s_0,u}(w),z,x,\rho) \rangle^2 dw\nu(dz) du }} .
\end{align*}
This shows the distribution equality after the same change of variable. $\square$

\bigskip

\bigskip

\bibliographystyle{abbrv}
\bibliography{biblio_stochsew}

\end{document}